\documentclass[11pt]{amsart}

\usepackage{amsfonts,amssymb,stmaryrd,amscd,amsmath,latexsym,amsbsy}

\usepackage{amssymb}
\usepackage{amsfonts}
\usepackage{latexsym}

\newtheorem{theorem}{Theorem}[section]
\newtheorem{lemma}[theorem]{Lemma}
\newtheorem{proposition}[theorem]{Proposition}
\newtheorem{corollary}[theorem]{Corollary}
\theoremstyle{definition}

\newcommand{\ben}{\begin{enumerate}}
\newcommand{\een}{\end{enumerate}}

\newcommand{\CC}{{\mathbb{C}}}
\newcommand{\ZZ}{{\mathbb{Z}}}

\theoremstyle{plain}

\newtheorem*{sol}{Solution}

\theoremstyle{definition}

\theoremstyle{remark}

\newcommand{\solu}[1]{\begin{sol}{\bf (\ref{#1})}}

\begin{document}

\title{Generalized double affine Hecke algebras of rank 1
and quantized Del Pezzo surfaces}

\author{Pavel Etingof}
\address{Department of Mathematics, Massachusetts Institute of Technology,
Cambridge, MA 02139, USA}
\email{etingof@math.mit.edu}

\author{Alexei Oblomkov}
\address{Department of Mathematics, Massachusetts Institute of Technology,
Cambridge, MA 02139, USA}
\email{oblomkov@math.mit.edu}

\author{Eric Rains}
\address{Department of Mathematics, University of California, Davis,
1 Shields Ave, Davis, CA 95616-8633, USA}
\email{rains@math.ucdavis.edu}

\maketitle

\section{Introduction}

Double affine Hecke algebras for reduced 
root systems were introduced by Cherednik 
\cite{Ch} in order to prove Macdonald conjectures. 
Double affine Hecke algebras of type 
$C^\vee C_n$ were introduced in the works of Noumi, Sahi, and Stokman
(\cite{NoSt,Sa,St}) as a generalization of Cherednik 
algebras of types $B_n$ and $C_n$, 
in order to prove Macdonald conjectures for Koornwinder
polynomials. 

The goal of this paper is to define and study new algebras $H(t,q)$, which
are generalizations of double affine Hecke algebras of type
$C^\vee C_n$ in the case $n=1$. To be more specific, 
fix a star-shaped simply laced affine Dynkin diagram $\widehat D$ 
(i.e., $\tilde D_4$, $\tilde E_6$, $\tilde E_7$, or
$\tilde E_8$). Let $m$ be the number of legs of $\widehat D$ and 
$d_j-1$, $j=1,...,m$, be the length of the $j$-th leg.
Then we define a family of algebras $H(t,q)$ 
depending on parameters $q\in \Bbb C^*$ and $t=(t_{kj})$,
$t_{kj}\in \Bbb C^*$, $k=1,...,m$, $j=1,...,d_k$, 
by generators 
$T_k$, $k=1,...,m$, with defining relations
\begin{equation}\label{defii1}
\prod_{j=1}^{d_k}(T_k-e^{2\pi ij/d_k}t_{kj})=0,\
k=1,...,m;\
\prod_{k=1}^mT_k=q.
\end{equation}
It follows from this definition 
that for $\widehat D=\tilde D_4$ we get exactly the double affine Hecke
algebra of type $C^\vee C_1$; on the other hand, 
in the case $\widehat D=\tilde E_{6,7,8}$ we get new algebras,
which are the main subject of this paper. 

It is obvious that if $t_{kj}=1$ and $q=1$, the algebra
$H(t,q)=H(1,1)$ is a group algebra of some group $G$. 
The group $G$ is defined by generators and relations, and 
is well known to be isomorphic to a 2-dimensional crystallographic 
group $\Bbb Z_\ell\ltimes \Bbb Z^2$, where $\ell=2,3,4,6$ in the
cases $\widehat D=\tilde D_4,\tilde E_6,\tilde E_7,\tilde E_8$,
respectively. Moreover, the algebra $H(1,q)$ is a twisted group
algebra of $G$. Thus, $H(t,q)$ is a deformation of the twisted
group algebra of $G$. We prove that if we regard 
$\log(t_{kj})$ as formal parameters then this deformation is flat 
({\it the formal PBW theorem}), and $H(t,qe^{\varepsilon})$ 
is the universal deformation of $H(1,q)$ if $q$ is not a root of
unity. We also prove a more delicate {\it algebraic PBW theorem},
which claims that (for numerical $t$, $q$) 
some filtrations on $H(t,q)$ have certain explicit
Poincar\'e series, independent of $t$ and $q$. 

It was shown by the second author (\cite{Ob}) that for 
$\widehat D=\tilde D_4$ and $q=1$ 
the algebra $H(t,q)$ is finite over its center $Z(t,q)$, and the 
spectrum of $Z(t,q)$ is an affine cubic surface, 
obtained from a projective one by removing three lines 
forming a triangle. Here we show that this result is 
valid also for $q$ being a root of unity, and generalize it
to the cases $\widehat D
=\tilde E_6,\tilde E_7,\tilde E_8$. In these cases, the spectrum of
$Z(t,q)$ ($q$ being a root of unity) turns out to be 
an affine surface $S(t,q)$ obtained from a projective del Pezzo 
surface $\overline{S(t,q)}$ 
of degrees $3,2,1$ respectively by removing a nodal $\Bbb
P^1$. This means that for $q\ne 1$, the spherical subalgebra 
$eH(t,q)e$ in the algebra $H(t,q)$ (where $e$ is the 
idempotent in $H(t,q)$ projecting to an eigenspace of the element
$T_3$ corresponding to the longest leg of $\widehat D$) 
should be viewed as an algebraic quantization of the surface $S(t,1)$
(with its unique up to scaling symplectic structure). 
Moreover, the algebraic PBW
theorem for $H(t,q)$ implies that the Rees algebra of $eH(t,q)e$
with respect to an appropriate filtration provides a 
quantization (in the sense of noncommutative algebraic geometry) 
of the projective Poisson surface $\overline{S(t,q)}$. 

The structure of the paper is as follows.  

In Section 2, we recall the basics about crystallographic
groups in the plane. 

In Section 3, we consider the 
twisted group algebra $B(q)=H(1,q)$ of a planar 
crystallographic group $G$, and deform it into an algebra
$\widehat{\bold H}(q)$, which is a version of $H(t,q)$ in which
$\log t_{kj}$ are formal parameters.  
We prove the formal 
PBW theorem for $\widehat{\bold H}(q)$, and formulate 
the results on the cohomology of $B(q)$ and on its 
universal deformation. 

In Section 4, we prove the results stated in Section 3. Namely, 
we compute the Hochschild cohomology of 
the twisted group algebra $B(q)$ and use the result to prove 
that $H(t,qe^{\varepsilon})$ is a universal deformation of
$B(q)$ if $q$ is not a root of unity. 

In Section 5, we define an increasing filtration on $H(t,q)$,
called the length filtration; its definition is based 
on a connection of $H(t,q)$ with a certain deformation 
of group algebras of affine Weyl groups (with lattice of rank 2). 
We show that the Poincar\'e series for this filtration is 
independent of $t,q$ (the algebraic PBW theorem). 
We then use this result to establish some general 
properties of $H(t,q)$, e.g. that the Gelfand-Kirillov dimension
of $H(t,q)$ is $2$. We also show that if $q$ is a root of unity
then the algebras $H(t,q)$ are PI. Finally, we show 
that if in addition $t$ is generic then $H(t,q)$ is an Azumaya
algebra, and the spectrum of the center $Z(t,q)$ of $H(t,q)$ 
is a smooth surface. 

To study finer structure of the algebras $H(t,q)$ at roots of
unity, one needs to define other filtrations on $H(t,q)$,
and prove the PBW theorem for them. 
This is done in Section 6. The proof of the PBW theorem 
is technical and relies on computer calculations; 
it is given in Section 8. We also compute 
the associated graded algebras attached to some of the
filtrations; the proof is again postponed till Section 8. 
In the second half of Section 6, we proceed to show that 
the spectrum of the center of $H(t,q)$ is an affine del
Pezzo surface. This shows that the spherical subalgebra 
$eH(t,q)e$ for $q\ne 1$ is a quantization 
of an affine del Pezzo surface, and also yields 
a linear algebra application given at the end of Section 6. 

In Section 7, using the Riemann-Hilbert correspondence, 
we define a homomorphism from 
a formal version of the generalized double affine Hecke algebra
to the completion of the deformed preprojective algebra 
of the quiver associated to the graph $\widehat D$.
This construction is similar to those used in \cite{CB2}.   
It allows us to define a holomorphic (but not algebraic) map from
the universal deformation of the Kleinian singularity 
$\Bbb C^2/\Gamma$ (where $\Gamma\subset SL_2(\Bbb C)$ is the
finite subgroup corresponding to the diagram $\widehat D$ via the McKay 
correspondence) to the family of surfaces $S(t,1)$. 
This is a local isomorphism of analytic varieties 
near $0\in \Bbb C^2/\Gamma$, which in the case 
$\widehat D=\tilde D_4$ encodes generic solutions of the 
Painlev\'e VI equation. 

In Section 8, we prove the results of Section 5, by writing 
presentations of $H(t,q)$ which are compatible with the filtrations.

Finally, in Section 9, we study more closely the 
surfaces $S(t,1)$. Namely, let $\bold G$ 
be the simple Lie group corresponding to the diagram $D$.
We show that the algebra $H(t,q)$ 
depends only on the projection of $t$ to the maximal torus of
$\bold T\subset \bold G$, and that the map $t\mapsto S(t,1)$ 
from $\bold T$ to the moduli space of affine del Pezzo surfaces 
is Galois and has Galois group isomorphic to the Weyl group $W$ 
of $\bold G$. Given the results of Section 7, 
this fact is in good agreement with the Arnold-Brieskorn 
theorem, saying that the monodromy group of a simple singularity 
is the Weyl group of the
corresponding Lie algebra. This also implies that the
coefficients of the equation of $S(t,1)$, as functions of $t$, 
 are polynomials of characters of irreducible representations of
$\bold G$, and we compute these polynomials explicitly.

We note that the computations of Sections 8 and 9 
are too complicated to be done by hand and were performed 
using a computer. More specifically, the third author wrote a
Magma code (\cite{Ma}) for computations in $H(t,q)$, 
which can be found at \cite{Ra}. 
We also remark that the results of Sections 1-5 and 7 are independent of
computer calculations. 

Finally, the paper contains two appendices. 
In Appendix 1, written by W. Crawley-Boevey and P. Shaw, 
it is shown that the algebra $H(t,q)$ 
is the ``spherical subalgebra'' (corresponding to the nodal
vertex idempotent) of the multiplicative preprojective algebra 
introduced in the paper \cite{CBS}. In Appendix 2, we use this
result to describe the structure of the multiplicative
preprojective algebras for affine starlike quivers. 
In particular, we show that if $q^\ell$ is a root of unity of
degree $N$ and $t$ is generic then the corresponding
multiplicative preprojective algebra is Azumaya of rank $hN$,
where $h$ is the Coxeter number of the corresponding Dynkin
diagram. 

{\bf Remark 1.} Some results of this paper can be extended to 
the case $n>1$, giving a generalization of double affine
Hecke algebras of type $C^\vee C_n$. This involves 
considering a flat 1-parameter deformation 
$H_n(t,q,k)$ of the algebra $\Bbb C[S_n]\ltimes H(t,q)^{\otimes
n}$. This is done in the subsequent paper \cite{EGO}; the
deformations considered there appear to provide quantizations 
of Hilbert schemes of Del Pezzo surfaces. 

{\bf Remark 2.} In \cite{GG}, Gan and Ginzburg define and study
rank $n$ analogs of preprojective algebras of quivers. It would be
interesting to define such analogs in the multiplicative
situation. In the case of affine quivers, they should 
have Gelfand-Kirillov dimension $2n$ and be finite over center
for special parameters. Also, for affine starlike quivers they
should be Morita equivalent to the algebras studied in \cite{EGO}
(see the previous remark).

{\bf Remark 3.} Finite dimensional representations
of the algebra $H(t,q)$ is essentially 
the same thing as solutions of the multiplicative
Deligne-Simpson problem, considered in \cite{CB2,CBS}. 
Thus the methods of \cite{CB2,CBS} can be used to obtain a 
classification of finite dimensional representations
of $H(t,q)$. For double affine Hecke algebras of type $A_1$, this
problem is solved in \cite{CO}, and for the (more general) type $C^\vee C_1$ in
\cite{OS}.  

{\bf Remark 4.} In \cite{VB}, M. Van den Bergh constructed
quantizations of del Pezzo surfaces with a (possibly singular) 
genus 1 curve removed, using the method of noncommutative
blowup; in the $E_6$ (=degree 3) case it was already 
done in \cite{LSV}. We expect that when the removed curve is 
a nodal $\Bbb P^1$, the algebras constructed in \cite{VB} 
are isomorphic to the spherical subalgebras $eH(t,q)e$
for $D=E_{9-d}$, where $d$ is the degree of the del Pezzo surface
($d=3,2,1$). Checking this should involve presenting both
algebras by generators and relations. 

{\bf Remark 5.} Let $S$ be a del Pezzo surface of degree $d\le 3$
with a genus 1 curve $E$ removed. Let $q$ be an automorphism of
$E$. In a forthcoming paper we plan to show that one can define an  
algebra $H_{S,E}(q)$ (depending continuously on $S,E,q$)
with an idempotent $e$ such that $eH_{S,E}(q)e$ 
is a quantization of $S$ (with its unique up to scaling
symplectic structure) and $H_{S,E}(q)=H(t,q)$ 
for a suitable $t$ when the curve $E$ is a nodal rational curve. 
Such algebras would provide elliptic deformations of 
generalized double affine Hecke algebras.  

{\bf Remark 6.} The algebras $H(t,q)$ are a special case of a much
more general class of algebras, which are flat deformations 
of group algebras of discrete groups, introduced in \cite{E} 
(for any group acting discretely on a complex manifold with
vanishing first and second Betti numbers), and in \cite{ER}
(for subgroups of even elements in Coxeter groups). These more
general deformations appear to be quite interesting (see e.g. \cite{EGO}), but 
are rather poorly understood at the moment. The authors 
plan to study them in subsequent papers. 

{\bf Acknowledgments.} The work of P.E. and A.O.
was partially supported by the NSF grant DMS-9988796
and the CRDF grant RM1-2545-MO-03. 
P.E. is very grateful to M. Artin for many useful explanations
about noncommutative algebraic geometry. We are also grateful 
to J. Starr for discussions about Del Pezzo surfaces, 
and to W. Crawley-Boevey, 
A. Malkin and M. Vybornov for explanations about preprojective
algebras of quivers. 

\section{Crystallographic groups in the plane and their twisted group
algebras}

Let $D$ be a simply laced Dynkin diagram, whose affinization $\widehat D$
has the structure of a star. That is, $\widehat D$ has a node
with $m$ legs growing out of it. Such diagrams $D$ are
$D_4$, $E_6$, $E_7$, and $E_8$ (the number $m$ is $4,3,3,3$, respectively).

Let $d_i$ be the length of the $i$-th leg of $\widehat D$ plus 1.
Consider the group $G$ generated by $T_i$, $i=1,...,m$,
with defining relations
$$
T_i^{d_i}=1, i=1,...,m,\ \prod_{i=1}^m T_i=1
$$
(fixing a cyclic ordering of legs).
Thus, for $D_4$ the group $G$ is generated by $a,b,c,d$ with
$$
a^2=b^2=c^2=d^2=1,\ abcd=1,
$$
for $E_6$ by $a,b,c$ with
$$
a^3=b^3=c^3=1,\ abc=1,
$$
for $E_7$ by $a,b,c$ with
$$
a^2=b^4=c^4=1,\ abc=1,
$$
for $E_8$ by $a,b,c$ with
$$
a^2=b^3=c^6=1,\ abc=1.
$$

It is well known (and easy to check) that $G$ is isomorphic
to the crystallographic group $\Bbb Z_\ell\ltimes \Bbb Z^2$,
where $\ell=2,3,4,6$ for $D_4,E_6,E_7,E_8$ respectively,
and the cyclic group $\Bbb Z_\ell$ acts on the lattice by rotations.

More specifically, we can view the group
$G$ as a group of affine transformations
of $\Bbb C$ using the following formulas for the action
of the generators.

In the $D_4$ case,
$$
a(z)=-z+1+{\rm i},\ b(z)=-z+1,\ c(z)=-z, \ d(z)=-z+{\rm i}.
$$

In the $E_6$ case,
$$
a(z)=\zeta(z+\zeta)-\zeta,\ b(z)=\zeta (z-1)+1, c(z)=\zeta z,
$$
where $\zeta=e^{2\pi {\rm i}/3}$.

In the $E_7$ case,
$$
a(z)=-z+1-{\rm i},\ b(z)={\rm i}(z-1)+1,\ c(z)={\rm i}z.
$$

In the $E_8$ case,
$$
a(z)=-z+1-\xi^2,\
b(z)=\xi^2(z-1)+1,\
c(z)=\xi z,
$$
where $\xi=e^{\pi {\rm i}/3}$.

Now let $q$ be an invertible variable, and
denote by $B$ the algebra generated
over $\Bbb C[q,q^{-1}]$
by $T_i$ with defining relations
$$
T_i^{d_i}=1, i=1,...,m,\ \prod_{i=1}^m T_i=q.
$$
If $q_0\in \Bbb C^*$, we can define
$B(q_0)=B/(q-q_0)$.

Let $q\in \Bbb C^*$ and $\widehat G$ be the central extension of
$G$ by $\Bbb Z$:
$\widehat G:=\Bbb Z_\ell\ltimes H$, where $H$ is the Heisenberg
group consisting of 3 by 3 upper triangular matrices with integer
entries and ones on the diagonal.  
Let $C$ be a generator of the center of $\widehat G$.
It is shown using the formulas above (see Section 4.1) that
$B(q)$ is a twisted group algebra
of $G$ -- the quotient of $\Bbb C[\widehat G]$ by the relation $C=q^\ell$.
In other words, $B(q)=\Bbb C[\Bbb Z_\ell]\ltimes A_{q^\ell}$,where
$A_Q$ is the $Q$-Weyl algebra generated by $X^{\pm 1}, P^{\pm 1}$
with the relation $PX=QXP$. Similar statements are true if $q$ is
a variable (i.e. if we work over $\Bbb C[q,q^{-1}]$). 

\section{Generalized double affine Hecke algebras
over formal series}

\subsection{Cohomology} We denote by $H^i(A,A)$ 
or simply $H^i(A)$ the Hochschild cohomology of an algebra
$A$. The Hochschild cohomology of an algebra $A$ with coefficients 
in a bimodule $M$ is denoted by $H^i(A,M)$. 

\begin{theorem} \label{coho}
Assume that $q\in \Bbb C^*$ is not a root of unity.
Then $H^0(B(q))=\Bbb C$,
$H^1(B(q))=0$, $H^2(B(q))=\Bbb C^{r+1}$,
$H^i(B(q))=0$ for $i>2$,
where $r$ is the rank of the Dynkin diagram $D$.
\end{theorem}

Let us now recall the definition of a {\it universal deformation}.
A flat $R$-algebra $A_R$ (with $R$ being a local commutative
Artinian algebra and $\mathfrak{m}\subset R$ the maximal ideal)
together with an isomorphism $A_R/\mathfrak{m}\simeq A$ is called
a flat deformation of $A$ over $S=Spec(R)$. A similar definition 
is made if $R$ is pro-Artinian. A flat deformation $A_R$  is a universal
deformation of $A$ if for every flat deformation $A_{\mathcal O(S)}$ of
$A$ over an Artinian base $S$ there exists a unique map $\tau$: $S\to
Spec(R)$ such that the isomorphism $A\simeq A_R/\mathfrak{m}$ lifts
to an isomorphism $A_{\mathcal O(S)}\simeq \tau^* A_R$.

It is well known that if $H^2(A,A)=E$ is a finite dimensional
vector space, and $H^3(A,A)=0$ then there exists a universal
deformation of $A$ parametrized by $E$ (i.e. with $R=\Bbb
C[[E]]$). Therefore, we have the following corollary. 

\begin{corollary}\label{coho1}
The universal deformation of $B(q)$ has $r+1$ parameters.
\end{corollary}

Now we will describe the universal deformation of $B(q)$ explicitly.

\subsection{Deformations of $B(q)$}

Let $t_{ij}$, $j=1,...,d_i$, $i=1,...,m$, be variables such that $\prod_j
t_{ij}=1$. Define $u_{kj}$ by the formula 
$u_{kj}=e^{2\pi j{\rm i}/d_k} t_{kj}$
for $k=1,...,m$.
We assume that $t_{ij}$ are formal, in the sense that
$t_{ij}=e^{\tau_{ij}}$, where $\tau_{ij}$ are formal
parameters. Let $t$ denote the collection of the variables
$t_{ij}$. Clearly, the number of independent variables among
them is $r$.

Define the algebra $\widehat {\bold H}$
to be (topologically) generated over $\Bbb C[q,q^{-1}][[\tau]]$ 
(where $\tau$ stands for the collection of variables $\tau_{ij}$)
by $T_k$, $k=1,...,m$, with defining relations
\begin{equation}\label{defii}
\prod_{j=1}^{d_k}(T_k-u_{kj})=0,\
k=1,...,m;\
\prod_{k=1}^mT_k=q.
\end{equation}
Sometimes we will use the notation $a,b,c,d$ for
$T_1,T_2,T_3,T_4$. 

This algebra of course depends on the Dynkin diagram $D$, but
in order to simplify notation we
will not write this dependence explicitly.
In the $D_4$ case, it is the double affine Hecke algebra
of type $C^\vee C_1$ of Sahi, Noumi and Stokman \cite{Sa,St,NoSt}.
So in the cases $E_6,E_7,E_8$, we get a generalization
of the double affine Hecke algebra. If $q_0\in \Bbb C^*$, we can also
define the algebra $\widehat{\bold H}(q_0):=\widehat{\bold
H}/(q-q_0)$ over $\Bbb C[[\tau]]$.

\begin{theorem}\label{pbw} (the formal PBW theorem)
The algebra $\widehat {\bold H}$ is a flat formal deformation
of $B$.
\end{theorem}

This immediately implies

\begin{corollary}\label{pbwq} For any $q\in \Bbb C^*$,
the algebra $\widehat {\bold H}(q)$ is a flat formal deformation
of $B(q)$.
\end{corollary}

We will show that if $q$ is not a root of unity then
this is the most general deformation. Namely, we have

\begin{theorem}\label{univ}
If $q$ is not a root of unity then
$\widehat {\bold H}(qe^{\varepsilon})$ (where $\varepsilon$
is a new formal parameter) is a universal deformation
of $B(q)$.
\end{theorem}

{\bf Remark.} In the $D_4$ case,
Theorem \ref{pbw} follows from the papers \cite{Sa,St,NoSt}.

\subsection{Proof of Theorem \ref{pbw}}

The proof is based on the following simple fact, which is often used 
for proving flatness of formal deformations.

\begin{lemma}\label{formdef} If $A_0$ is an
algebra over $\Bbb C$, $A$ a formal deformation of $A_0$ (over $\Bbb
C[[t_1,...,t_n]]$), and $M_0$ is a faithful $A_0$-module which
can be flatly deformed to an $A$-module, then $A$ is a {\bf flat}
formal deformation of $A_0$. 
\end{lemma}

For every element $g\in G$, fix its presentation as a product of
$T_i$, and denote by $b_g^0$ the same product in
the algebra $B$. Then $\lbrace{b_g^0\rbrace}$ is a basis of $B$
labeled by group elements $g\in G$ 
(it is independent of the choice of the presentations 
up to scaling by powers of $q$). Let $b_g$ be some lifts of $b_g^0$ to
$\widehat{\bold H}$. Let $J$ be the maximal ideal in $\Bbb
C[[\tau]]$.

Assume that $\widehat{\bold H}$ is not flat.
Then there exists $n>0$ and $g_1,...,g_k\in G$,
$\alpha_1,...,\alpha_k\in \Bbb C[q,q^{-1}][[\tau]]/J^n$
(not all zero) such that
$\sum \alpha_jb_{g_j}=0$ in $\widehat{\bold H}/J^n$.
This relation is nontrivial if we reduce it
modulo $q-q_0$ for all but finitely many $q_0$.
Hence $\widehat{\bold H}(q)$ is not a flat deformation
of $B(q)$ for all but finitely many $q$.

Thus it is sufficient to establish that $\widehat{\bold H}(q)$ is
flat in the case when
$q$ is a root of unity.
We will assume that $q^\ell$ is a root of unity of degree $N$, where $N$
is a positive integer.

In this case, the algebra $B(q)$ is a semidirect product
of $\Bbb Z_\ell$ with the q-Weyl algebra (or quantum torus) $A_{q^\ell}$. 
The algebra $A_{q^\ell}$ is well known to be an Azumaya algebra of rank $N$. 
In particular, the center of $B(q)$ is $Z=\Bbb C[T/\Bbb Z_\ell]$,
and $B(q)\otimes_Z \overline Q_Z={\rm Mat}_{\ell N}(\overline Q_Z)$
(where $Q_Z$ is the field of fractions of $Z$, and $\overline
Q_Z$ is the algebraic closure of $Q_Z$).
So $B(q)$ admits a 2-parameter family of
irreducible $\ell N$-dimensional representations, whose direct
sum is faithful. By Lemma \ref{formdef}, 
to prove our theorem, it is sufficient to show that
these representations can be deformed to representations of
$\widehat{\bold H}(q)$.

Let $V$ be an irreducible $\ell N$-dimensional matrix representation
of $B(q)$. Let ${\bold C}_k$,
$k=1,...,m$, be the conjugacy class
of ${\rm diag}(u_{kj})\otimes {\rm Id}_{N\ell/d_k}$ in $GL_{\ell N}$;
it is a smooth algebraic variety defined over $\Bbb C[[\tau]]$.
Consider the scheme ${\bold Y}$ of m-tuples
$(T_1,...,T_m)$ lying in the formal neighborhood of the orbit of
$V$ (under changes of basis), 
such that $T_k\in {\bold C}_k$ and $\prod_k T_k=q {\rm Id}_{\ell
N}$. Our job 
is to show that the structure ring $O_{\bold Y}$ is flat over
$\Bbb C[[\tau]]$.

Let $Y$ be the reduction of ${\bold Y}$ modulo the maximal
ideal in $\Bbb C[[\tau]]$. By Schur's lemma,
$Y$ admits a free action of $PGL_{\ell N}(\Bbb C)$, and the
quotient is a 2-dimensional formal polydisk. Thus
$Y$ is smooth and has dimension $\ell^2 N^2+1$.

On the other hand, let us compute the ``expected dimension'' of
${\bold Y}$, i.e. the dimension of the ambient space minus the
number of equations.
Fixing a matrix $T_k\in {\bold C}_k$ amounts to fixing
$d_k$ subspaces in an $\ell N$-dimensional linear space
of dimension $\ell N/d_k$ which add up to the whole space.
Thus, ${\bold C}_k$
has dimension
$$
D_k=\ell^2 N^2(1-\frac{1}{d_k}).
$$
On the other hand, the number of equations in the condition
$\prod_k T_k= q {\rm Id}_{\ell N}$ is $\ell^2 N^2-1$ (since
the determinant of the product is fixed). Thus the expected
dimension is
$$
{\Bbb D}=\sum_k D_k-\ell^2
N^2+1=\ell^2N^2\sum_{k=1}^m(1-\frac{1}{d_k})-\ell^2 N^2+1
$$

But $\sum_{k=1}^m (1-\frac{1}{d_k})=2$ (as $\widehat D$ is an
affine diagram). Thus,
${\Bbb D}=\ell^2 N^2+1$. The expected dimension of $Y$ is
obviously the same.

Thus, the expected dimension of $Y$ coincides with its actual
dimension. This implies that $Y$ is a complete intersection, 
and therefore so is ${\bold Y}$. Since ${\bold Y}$ is obtained
from $Y$ by deforming its equations, and $Y$ is a complete
intersection, we conclude that $\bold Y$ is a flat deformation of
$Y$ (in fact, it is, moreover, a trivial deformation, since $Y$ is smooth).
The theorem is proved.

\section{Proofs of Theorems \ref{coho},\ref{univ}}

\subsection{Homology and cohomology of $B(q)$}
In this section we prove Theorems~\ref{coho},\ref{univ}. We will use arguments
similar to the arguments from \cite{Ob}.

Let us describe  the
isomorphism between the semidirect product $B=\Bbb C[\Bbb
Z_\ell]\ltimes A_{q^\ell}$ and the
algebra $B(q)$. For 
brevity we use the symbol $D_q$ for the algebra
$A_{q^\ell}$.

In the case $\ell=2$ in the formulas defining $B(q)$ we
have $d_1=d_2=d_3=d_4=2$ and we can choose the alternative set of
generators $P,X,s$:
$$
 X=T_1T_2,\quad P=T_2T_3,\quad s=T_4.$$
These elements generate $B(q)$ modulo the relations:
$$ PX=q^2 XP,\quad s^{-1}Ps=P^{-1},\quad s^{-1} X s=X^{-1},
\quad s^2=1.$$

In the case $\ell=3$ in the definition of $B(q)$ we have
$d_1=d_2=d_3=3$, and the alternative system of generators is:
\begin{equation}\label{gensforl=3}
 X=T_3T_1^{-1}, \quad P=T_1 T_2^{-1}, s=T_3. 
\end{equation}
These elements generate $B(q)$ modulo the relations:
\begin{equation}\label{relsforl=3}
 PX=q^3 XP,\quad s^{-1}Xs=q^{-1} P^{-1}X^{-1},\quad
s^{-1} P s=q^2 X,\quad s^3=1.
\end{equation}

In the case $\ell=4$ in the definition of $B(q)$ we have
$d_1=2,d_2=d_3=4$, and the alternative system of generators is:
\begin{equation}\label{gensforl=4}
 X=T_2^2 T_1,\quad  P=T_3^2 T_1, \quad s=T_3.
\end{equation}
These elements generate $B(q)$ modulo the relations:
\begin{equation}\label{relsforl=4}
PX=q^4 XP,\quad s^{-1}Xs=q^{-2}P^{-1},\quad s^{-1}Ps=q^{2}X,\quad
s^4=1.
\end{equation}

In the case $\ell=6$ in the definition of $\Bbb C_q[G]$ we have
$d_1=2, d_2=3, d_3=6$, and the alternative system of generators
is:
\begin{equation}\label{gensforl=6}
X=T_3^3T_1,\quad P=T_3^{-2}T_2, \quad s=T_3.
\end{equation}
These elements generate $\Bbb C_q[G]$ modulo the relations:
\begin{equation}\label{relsforl=6}
PX=q^6XP, \quad s^{-1}Xs=q^{-2}XP^{-1},\quad s^{-1} Ps=qX,\quad
s^6=1.
\end{equation}

\begin{theorem}\label{HH} If $q$ is not a root of unity then
for $\ell=2,3,4,6$ we have
\begin{gather*}
H^2(\mathbb{Z}_\ell\ltimes D_q)=H_0(\mathbb{Z}_\ell\ltimes D_q)=
\mathbb{C}^{r+1},\\
 H^1(\mathbb{Z}_\ell\ltimes D_q)=H_1(\mathbb{Z}_\ell\ltimes D_q)=0,\\
H^0(\mathbb{Z}_\ell\ltimes D_q)=
H_2(\mathbb{Z}_\ell\ltimes D_q)=\mathbb{C},\\
H_{>2}(\mathbb{Z}_\ell\ltimes D_q)=H^{>2}(\mathbb{Z}_\ell\ltimes D_q)=0.
\end{gather*} with $r=4,6,7,8$ respectively.
\end{theorem}

We  prove this theorem using the technique from \cite{AFLS}.

Recall (\cite{EtOb}) that $D_q\in VB(2)$, i.e., 
there exists an isomorphism of bimodules
$\zeta: H^2(D_q,D_q\otimes D_q^{opp})\to D_q$,
 where $D_q^{opp}$ is the algebra $D_q$ with the opposite multiplication.

\begin{lemma}\label{unim} The isomorphism $\zeta$ is 
$\Bbb Z_\ell$-equivariant. 
\end{lemma}

\begin{proof} Since the center of $D_q$ is trivial, 
an isomorphism $\zeta$ is unique up to scaling, hence 
$\Bbb Z_\ell$ must act on $\zeta$ by a character
$\chi: G\to \Bbb C^*$. Then by Van den Bergh's theorem 
(\cite{VB1,VB2}), $H_0(D_q,D_q)=H^2(D_q,D_q)\otimes \chi$ as 
$\Bbb Z_\ell$-modules. But $H_0(D_q,D_q)$ 
is the trivial $\Bbb Z_\ell$-module, since 
$D_q$ has a unique trace (up to scaling), sending $X^iP^j$ to
$\delta_{i0}\delta_{j0}$, and this trace is clearly fixed under
$\Bbb Z_\ell$. On the other hand, $H^2(D_q,D_q)$ is 
1-dimensional (as is easily seen from the Koszul resolution, see
below), and spanned by the class defined by the deformation of $D_q$ into
$D_{qe^\varepsilon}$. Since $\Bbb Z_\ell$ acts on
$D_{qe^{\varepsilon}}$, we see that this class is also fixed 
under $\Bbb Z_\ell$. Thus, $\chi=1$ and the lemma is proved.    
\end{proof} 

Lemma \ref{unim} and Proposition 3.5 from \cite{EtOb}
imply that there is an isomorphism
between the Hochschild homology
$H_i(\mathbb{Z}_\ell\ltimes D_q)$ and Hochschild cohomology
$H^{2-i}(\mathbb{Z}_\ell\ltimes D_q)$.
Thus it suffices to calculate the Hochschild homology
$H_*(\mathbb{Z}_\ell\ltimes D_q)$.

\subsection{The decomposition of the Hochschild homology}\label{spseq}

There is a natural structure of a $\ZZ_{\ell}$-module on the
homology $H_i(D_q,gD_q )$, where $g\in
\ZZ_\ell$. More precisely, there is an action of $\ZZ_\ell$ on the
standard Hochschild complex for $H_i(D_q,gD_q)$ by the formulas:
$$g\cdot(m\otimes a_1\otimes\dots \otimes a_r)=m^{g}\otimes
a_1^g\otimes\dots\otimes a_r^g.$$

 Proposition 3.1 from the paper \cite{AFLS} implies:
\begin{proposition}There is a decomposition:
$$H_*(\mathbb{Z}_\ell\ltimes D_q)=\oplus_{g\in\mathbb{Z}_\ell}
 H_*(D_q, gD_q)^{\mathbb{Z}_\ell}.$$
\end{proposition}
For calculation of $H_*(D_q,gD_q)$ we will use the Koszul
resolution.

\subsection{Calculation of $H_*(D_q,gD_q)$}
 Let us denote by $D^e_q$ the algebra $D_q\otimes D_q^{opp}$.
The elements $p=P\otimes P^{-1}-1$, $x=X\otimes X^{-1}-1$ commute
and $D_q^e/{\bf I}=D_q$, where ${\bf I}=(x,p)$ is the
$D_q^e$-submodule generated by these elements. Hence the
corresponding Koszul complex yields a free resolution $W_*$ of
the $D^e_q$-module $D_q$:
$$D^e_q\stackrel{{\bold d}_1}{\to} D^e_q\oplus
D^e_q\stackrel{{\bold d}_0}{\to}D^e_q\stackrel{\mu}{\to} D_q,$$ where
$\mu(X^iP^j\otimes P^{j'}X^{i'})=X^i P^{j+j'}X^{i'}$ and for
$z=z_1\otimes z_2$ we have ${\bold d}_0(z,0)=zp=z_1 P\otimes P^{-1}
z_2-z_1\otimes z_2$, ${\bold d}_0(0,z)=zx=z_1X\otimes X^{-1}z_2-z_1\otimes
z_2$, ${\bold d}_1(z)= (zx,-zp)=(z_1X\otimes X^{-1}z_2-z_1\otimes
z_2,-z_1P\otimes P^{-1} z_2+z_1\otimes z_2).$

Using the Koszul complex for $D_q$ we prove the following
\begin{lemma}\label{HH0} Suppose that we have an automorphism $g$ 
of the algebra $D_q$ given by the formulas:
\begin{gather*}
X^g=q^{b_1}X^{g_{11}} P^{g_{21}},\\
P^g=q^{b_2}X^{g_{12}}P^{g_{22}},
\end{gather*}
(where $det(g_{ij})=1$), and suppose that the map $g-1$: $\mathbb{Z}^2\to \mathbb{Z}^2$ is
injective. Let $v^1,\dots,v^k\in \mathbb{Z}^2$ be vectors such
that $v^i-v^j\notin (g-1)\mathbb{Z}^2$ for $i\ne j$ and
$k=|\mathbb{Z}^2/(g-1)\mathbb{Z}^2|$. Then we have:
\begin{equation*}
H_{>0}(D_q,gD_q)=0,
\end{equation*}
and
\begin{equation*}
H_0(D_q,gD_q)=gD_q/[D_q,gD_q]=\oplus_{i=1}^k \mathbb{C}\langle
X^{v^i_1}P^{v^i_2}+[D_q,gD_q]\rangle.
\end{equation*}
\end{lemma}
\begin{proof}
Identifying the vector spaces $gD_q$ and $D_q$ by sending $gx\in gD_q$
to  $x\in D_q$ and taking the tensor product of the resolution $W_*$
and $gD_q$ over $D^e_q$, we get the complex
$$ D_q\stackrel{\hat{{\bold d}_1}}{\to} D_q\oplus
D_q\stackrel{\hat{{\bold d}_0}}{\to} D_q,$$ where $\hat{{\bold d}_1}(z)=
(X^gzX^{-1}-z,-P^gzP^{-1}+z)$, $\hat{{\bold d}_0}(w,z)=P^gwP^{-1}-w+
X^gzX^{-1}-z.$ The homology of this complex is exactly the homology
$H_*(D_q,gD_q)$.

It is easy to check that $\hat {\bold d}_1$ is injective, so
we have $H_2(D_q,gD_q)=0$.
So it remains to compute $H_1$ and $H_0$.

Let us write the maps $\hat{{\bold d}_0}$, $\hat{{\bold d}_1}$ in terms of the PBW
basis in $D_q$.  A direct calculation shows that
\begin{gather*}
\hat{{\bold d}_0}(w,y)=\sum_{i,j}
(\delta^{(1)}c_2+\delta^{(2)}c_1)(i,j)X^iP^j,\\
\hat{{\bold d}_1}(z)=(\sum_{i,j} \delta^{(1)}c(i,j)X^iP^j,
-\sum_{i,j}\delta^{(2)}c(i,j) X^iP^j),
\end{gather*}
where $w=\sum_{i,j} c_1(i,j)X^iP^j$, $y=\sum_{i,j}c_2(i,j)X^i
P^j$, $z=\sum_{i,j} c(i,j) X^iP^j$ and
\begin{gather*}
(\delta^{(1)}c)(i,j)=exp(h(g_{21}i-j+(1-g_{11})g_{21}+b_1))c(
(i,j)-w_1)-c(i,j),\\
(\delta^{(2)}c)(i,j)=exp(h(g_{22}i-g_{22}g_{12}+b_2))c((i,j)-w_2)-
c(i,j),
\end{gather*}
with $w_1=(g_{11}-1,g_{21})$, $w_2=(g_{12},g_{22}-1)$ being a
basis of the lattice $(g-1)\mathbb{Z}^2$ and $q=exp(h)$.

The operations $\delta^{(i)}$, $i=1,2$  preserve the space
$F_{fin}$ of the functions on $\mathbb{Z}^2$ with finite support,
and they obviously commute. These operations are discrete
analogues of partial differentiations, and the image of $\delta^{(i)}$
could be described in terms of the discrete analog of integration:
\begin{gather*}(I^{(1)}c)(i,j)=\sum_{(m,n)\in
(i,j)+kw_1,k\in\mathbb{Z}} exp(hs(m,n))c(m,n),\\
(I^{(2)}c)(i,j)=\sum_{(m,n)\in (i,j)+k w_2,k\in\mathbb{Z}}
exp(hs(m,n))c(m,n),
\end{gather*}
where $s$: $\mathbb{Z}^2\to \mathbb{Z}$ is any function such that
$I^{(i)}$ satisfy the equations
$$ I^{(1)}\delta^{(1)}c=0,\quad  I^{(2)}\delta^{(2)}c=0,$$
for any function $c\in F_{fin}$. The last equations are equivalent
to the system:
\begin{gather*}
s((i,j)+w_1)=-i g_{21}+g_{21}+j+s(i,j)-b_1,\\
s((i,j)+w_2)=s(i,j)-ig_{22}-b_2.
\end{gather*}
This system has a $k$-dimensional (affine) space of solutions. Indeed, if we
fix the values of $s(v_i)$, $i=1,\dots,k$ then the value of $s$ at
the point $v_i+mw_1+nw_2$ could be found from the system. In
particular the solution of the system normalized by the
condition $s(v_i)=0$ has the form
\begin{gather*}s(v)=s_i(m,n),\mbox{ for } v=v_i+mw_1+nw_2,
 \end{gather*} where
% exp(hs(v))=0 \mbox{ if } v-v_i\notin (g-1)\mathbb{Z}^2,
\begin{gather*}
s_i(m,n)=\frac{(2-g_{11})g_{21}m^2}{2}+g_{22}(1-g_{11})mn-
\frac{g_{11}g_{22}n^2}{2}+\\
(\frac{g_{12}g_{22}}{2}-b_2-g_{22}v_i^1)n+(g_{21}(g_{11}-v_i^1)+
v_i^2-b_1)m.
\end{gather*}

Thus we have the following description of the image of
$\delta^{(i)}$: $$ {\rm Im}\delta^{(i)}={\rm Ker} I^{(i)}.$$ Having this
description we can show that ${\rm Ker} \hat{\bold d}_0\subset
{\rm Im}\hat{\bold d}_1$.
Indeed if $(w,y)\in {\rm Ker} \hat{\bold d}_0$, $w=\sum_{i,j} c_1(i,j) X^i
P^j$, $y=\sum_{i,j}c_2(i,j)X^iP^j$ then
$\delta^{(1)}c_2=-\delta^{(2)}c_1$. As $I^{(1)}$ commutes with
$\delta^{(2)}$, we have:
\begin{gather*}
0=I^{(1)}\delta^{(1)}c_2=-I^{(1)}\delta^{(2)}c_1=-\delta^{(2)}I^{(1)}c_1.
\end{gather*}

As $c_1\in F_{fin}$, the  equation implies $I^{(1)}c_1=0$, hence
$c_1=\delta^{(1)} c$ for some $c\in F_{fin}$. A similar
calculation shows that  $c_2=\delta^{(2)} c'$ for some $c'\in
F_{fin}$. Moreover, the equation:
\begin{equation*}
\delta^{(1)}\delta^{(2)}(c+c')=\delta^{(2)}c_1+\delta^{(1)}c_2=0,
\end{equation*}
implies $c+c'=0$. Thus $H_1(D_q,gD_q)=0$.

Let us now prove that $H_0(D_q,gD_q)=\mathbb{C}^k$. Indeed, it is easy
to see that ${\rm Im}\hat{{\bold d}_0}=\cap_{s=1}^k{\rm Ker} I_s$, where $I_i$:
$F_{fin}\to\mathbb{C}$:
$$ I_i=\sum_{(m,n)\in\mathbb{Z}^2}q^{s(m,n)} f(v_i+mw_1+nw_2).$$
Thus $\dim H_0(D_q,gD_q)\le k$.

On the other hand, the
vectors $gX^{v_1^s}P^{v_2^s}$, 
$s=1,\dots,k$ are linearly independent modulo the
subspace $[D_q,gD_q]\subset gD_q$ because the subspace
$[D_q,gD_q]$ is spanned by the vectors of the form
$g(X^iP^j-q^{f(i,j,u)}X^{i+u_1}P^{j+u_2})$ where $u\in
(g-1)\mathbb{Z}^2$ and $f(i,j,u)$ is some function.
Thus $H_0(D_q,gD_q)=gD_q/[D_q,gD_q]=\mathbb{C}^k$, and the lemma
is proved.
\end{proof}

\begin{corollary} If $q$ is not a root unity and
$g\in \mathbb{Z}_\ell$ is not the unit element then $H_{>0}(D_q,
gD_q)=0$ and
\begin{gather*}
H_0(D_q,sD_q)=\mathbb{C}^4, \mbox{ for }\ell=2,\\
H_0(D_q,sD_q)=H_0(D_q,s^2D_q)=\mathbb{C}^3, \mbox{ for }\ell=3,\\
H_0(D_q,sD_q)=H_0(D_q,s^3D_q)=\mathbb{C}^2,\quad
H_0(D_q,s^2D_q)=\mathbb{C}^4, \mbox{ for }\ell=4,\\
H_0(D_q,sD_q)=H_0(D_q,s^5D_q)=\mathbb{C},\quad
H_0(D_q,s^2D_q)=H_0(D_q,s^4D_q)=\mathbb{C}^3, \\
H_0(D_q,s^3D_q)=\mathbb{C}^4,\mbox{ for } \ell=6,
\end{gather*}
where $s$ is a generator of $\Bbb Z_\ell$. 
\end{corollary}

\subsection{Calculation of $H_*(D_q,D_q)^{\mathbb{Z}_\ell}$}
Using the Koszul resolution we can easily calculate the homology of
$D_q$ (see for example section 4 of \cite{Ob}).

\begin{lemma} If $q$ is not a root of unity then
$$H_0(D_q,D_q)=H_2(D_q,D_q)=\mathbb{C},\quad
H_1(D_q,D_q)=\mathbb{C}^2.$$ Moreover, for all $\ell$ we have
$$H_0(D_q,D_q)^{\mathbb{Z}_\ell}=H_2(D_q,D_q)^{\mathbb{Z}_\ell}
=\mathbb{C},\quad H_1(D_q,D_q)^{\mathbb{Z}_\ell}=0.$$
\end{lemma}

\begin{proof}
Indeed the first statement follows from a simple
calculation with Koszul resolution; this calculation is done for
example in section 4 of \cite{Ob}.

Let us now prove the second statement.
Recall that by Lemma \ref{unim} and the results of \cite{VB1,VB2}
the space $H_1(D_q,D_q)$ is $\Bbb Z_\ell$ equivariantly
isomorphic  to the space
$$H^1(D_q,D_q)={\rm Der}(D_q)/\langle [x,\cdot],x\in D_q \rangle,$$
where ${\rm Der}(D_q)$ is the space of the derivations of $D_q$
i.e. the space of the $\mathbb{C}$-linear maps $d: D_q\to D_q$
with the property $d(xy)=d(x)y+xd(y)$. From the last description
we see that $H^1(D_q,D_q)$ is spanned\footnote{Note that this is only true if $q$
is not a root of unity; if $q$ is a root of unity then the spaces
$H^i(D_q,D_q)$ are infinite dimensional for $i=0,1,2$} by two derivations ${\bold d}_1$,
${\bold d}_2$: ${\bold d}_1(X^i P^j)=iX^iP^j$, ${\bold
d}_2(P^jX^i)=jP^jX^i$. The action
of the group $\mathbb{Z}_\ell$ on the space ${\rm Der}(D_q)$ is given
by the formula $(g\cdot d)(x)=g^{-1}(d(g(x)))$. Thus the action
of element $g\in\mathbb{Z}_\ell$, $\ell=2,4$ of order $2$ is given
by the formula:
$$g({\bold d}_1)=-{\bold d}_1,\quad g({\bold d}_2)=-{\bold d}_2,$$
hence $H^1(D_q,D_q)^{\mathbb{Z}_4}\subset
H^1(D_q,D_q)^{\mathbb{Z}_2}=0$. The action of the element
$g\in\mathbb{Z}_\ell$, $\ell=3,6$ of order $3$ is given by the
formula:
$$ g({\bold d}_1)=-{\bold d}_1+{\bold d}_2,\quad g({\bold d}_2)=-{\bold d}_1,$$
hence $H^1(D_q,D_q)^{\mathbb{Z}_6}\subset
H^1(D_q,D_q)^{\mathbb{Z}_3}=0$.

As $H^0(D_q,D_q)=Z(D_q)=\mathbb{C}$ is isomorphic to
$H_2(D_q,D_q)=\mathbb{C}$ and the action of $\mathbb{Z}_\ell$ on
the center is trivial, we get
$H_2(D_q,D_q)^{\mathbb{Z}_\ell}=\mathbb{C}$. Analogously,
we have seen in the proof of Lemma \ref{unim} 
that $H^2(D_q,D_q)$ and $H_0(D_q,D_q)=D_q/[D_q,D_q]$ 
are trivial $\Bbb Z_\ell$-modules. 
\end{proof}

\subsection{Proof of theorem~\ref{HH}}
 To complete the proof we only need to calculate the action of
$\mathbb{Z}_\ell$ on the cohomology $H_0(D_q,s^iD_q)$,
$i=1,\dots,\ell-1$. Let us recall that the action of the group
$\mathbb{Z}_\ell$ on $H_0(D_q,s^iD_q)=s^iD_q/[s^iD_q,D_q]$ is
induced by the conjugation action of $\mathbb{Z}_\ell$ on
$s^iD_q$. As $s^i(x)-x\in [s^iD_q,D_q]$ we get that the action of $s^i$
on $H_0(D_q,s^iD_q)$ is trivial and we only need to calculate
$H_0(D_q,s^2D_q)^{\mathbb{Z}_4}$ in the case $\ell=4$ and
$H_0(D_q,s^2D_q)^{\mathbb{Z}_6}$,
$H_0(D_q,s^3D_q)^{\mathbb{Z}_6}$,
$H_0(D_q,s^4D_q)^{\mathbb{Z}_6}$. Let us prove
$H_0(D_q,s^2D_q)^{\mathbb{Z}_4}=\mathbb{C}^3$. Indeed, this follows
by calculation from the fact that we have the equalities: $$ s\cdot
s=s,\quad s\cdot (s^2X)=q^{-2}s^2P, \quad s\cdot(s^2 P)=q^2 s^2X,
\quad s\cdot(s^2XP)=s^2XP$$ modulo $[D_q,s^2D_q].$ An analogous
calculation shows that
$$H_0(D_q,s^2D_q)^{\mathbb{Z}_6}=H_0(D_q,s^4D_q)^{\mathbb{Z}_6}=
\mathbb{C}^2,\quad H_0(D_q,s^3D_q)^{\mathbb{Z}_6}=\mathbb{C}^2.$$
Thus Proposition~\ref{spseq} implies the theorem.

\subsection{Infinitesimal deformations}
In this subsection we prove Theorem~\ref{univ}.

Let $q\in \Bbb C^*$, and let $({\tau},h)$, 
where $\tau=(\tau_{kj})\in \Bbb C^r$ and $h\in \Bbb C$, 
be a nonzero vector. Let
$H'$ be the algebra over the ring $\Bbb
C[\varepsilon]/\varepsilon^2$
of dual numbers, generated 
by $T_i$ with defining relations 
(\ref{defii}), with $t_{ij}=e^{\varepsilon \tau_{ij}}$, 
and $q$ replaced by $qe^{h\varepsilon}$. 
Theorem \ref{univ}
follows from Theorem~\ref{HH} and the following
lemma (in which $q$ is allowed to be a root of unity).

\begin{lemma}\label{nontriv}
 There is no isomorphism of $\CC[\varepsilon]/(\varepsilon^2)$-algebras
between $H'$ and
$B(q)\otimes_{\CC}\CC[\varepsilon]/(\varepsilon^2)$ which
is  equal to the identity map modulo the ideal $(\varepsilon)$.
\end{lemma}

In the case $\ell=2$ this lemma was proved in \cite{Ob}. We show
how we can modify this proof for the cases $\ell=3,4,6$.

In the proof of this lemma we use the following description of the
algebra $H'$: it is generated by $T_i$ with defining relations
\begin{gather*}
T_i^{d_i}=1+\varepsilon\sum_{k=1}^{d_i-1}L^i_k({\tau})
T_i^{d_i-k},\quad
i=1,2,3\\
T_1T_2T_3=qe^{h\varepsilon},
\end{gather*}
where
$L_k^i({\tau})$
are appropriate formal series in ${\tau}$.

\begin{proof}[Proof of the Lemma~\ref{nontriv}]

First we explain the proof in the case $\ell=3$. Let us denote by
$\phi$ the natural isomorphism of vector spaces
$$B(q)\to
\varepsilon H'.$$
Assume that there is an isomorphism between $H'$ and
$B(q)\otimes_{\CC}\CC[\varepsilon]/(\varepsilon^2)$
lifting the identity. A direct calculation shows that the
following equation holds in $H'$:
\begin{multline}\label{modeps}
PX-q^3XP=\varepsilon(3hq^3XP-s(q^{-2}L^1_2X+q^2L^2_2P^{-1}X+L^3_3PX)+\\
s^2(q^2L^1_1+qL^2_1X+L^3_1PX)),
\end{multline}
where $P,X$ are given by formulas (\ref{gensforl=3}).
 Hence there exist elements $X',P'\in B(q)$ 
such that the elements $P+\varepsilon P'$
and $X+\varepsilon X'$ of $B(q)[\varepsilon]/\varepsilon^2$
satisfy  the equation (\ref{modeps}) modulo
$\varepsilon^2$. If we write the elements $X'$, $P'$ in the PBW
basis:
\begin{equation*}
X'=\sum_{i=0}^2 s^i\Phi^X_i,\quad P'=\sum_{i=0}^2 s^i\Phi^P_i,
\end{equation*}
with $\Phi^X_i,\Phi^P_i\in\mathbb{C}_q[X^{\pm 1},P^{\pm 1}]$, then
the equality of the 
coefficients before $\varepsilon$ in the equation (\ref{modeps})
yields:
\begin{gather}
\Phi_0^PX-q^3X\Phi_0^P=3hq^3XP,\label{eq0}\\
\Phi_1^P X-q^3 X^s\Phi_1^P+P^s\Phi_1^X-q^3\Phi_1^XP=
-(q^{-2}L^1_2X+q^2L^2_2P^{-1}X+L^3_2PX),\label{eq1}\\
\Phi_2^P-q^3X^{s^2}\Phi_2^P+P^{s^2}\Phi^X_2-q^3\Phi_2^XP=
q^2L^1_1+q^2L^2_1X+L^3_1PX.\label{eq2}
\end{gather}
Here $X^{s^j}:=s^{-j}Xs^j$ and 
$P^{s^j}:=s^{-j}Ps^j$.

If we expand the LHS of (\ref{eq0}) in terms of the PBW basis $X^iP^j$, we
see that the coefficient before $XP$ is zero, hence $h=0$.
Similarly to the description of ${\rm Im}\hat{\bold d}_0$ from the proof of
Theorem~\ref{HH} we see that the coefficients $c_{ij}$ of
expansion $\sum c_{ij}X^iP^j$ of LHS of (\ref{eq1})  satisfy the
equations $$I_n(c)=\sum_{(i,j)\in v_n+(s-1)\mathbb{Z}^2}
q^{f_n(i,j)}c_{i,j}=0,$$ $n=1,2,3$ where $v_i\in\mathbb{Z}^2$ are
distinct modulo $(s-1)\mathbb{Z}^2$ and $f_n$ is a quadratic
expression of $i,j$. At the same time, we see that the RHS of
(\ref{eq1}) has the form $\sum_{n=1}^3c_n(q,\tau)X^{w^1_n}P^{w^2_n}$,
and vectors $w_i$, $i=1,2,3$ are distinct modulo
$(s-1)\mathbb{Z}^2$. Thus we proved that $L^1_2=L^2_2=L^3_2=0$ for
all $i$. Analogously, considering equation (\ref{eq2}) we prove
that $L^1_1=L^2_1=L^3_1=0$.

Thus  we get that the vectors
$\tau_i=(\tau_{i1},\tau_{i2},\tau_{i3})$ satisfy the equation:
$$ L_j^i(\tau_i)=0,\quad j=1,2,3.$$
It is easy to see that Jacobi matrix of this system is
nondegenerate 
at zero, so we have $\tau_i=0$, $i=1,2,3$.

The proof in the case $\ell=4,6$ is the same, except that we need to use
instead of formula (\ref{modeps}) the formulas:
\begin{multline*}
PX-q^4XP=\varepsilon(4hq^4XP+s(q^{m_{11}}L^3_3-q^{m_{12}}L^2_3P)+s^2
(q^{m_{21}}L^2_2P^2-q^{m_{22}}L_2^3P+\\
L_1^1(q^{m_{23}}P-
q^{m_{24}}X^{-1}P^2))+s^3(q^{m_{31}}L_1^3P^2-q^{m_{32}}L_1^2X^{-1}P^2)),
\end{multline*}
for the case $\ell=4$, and
\begin{multline*}
PX-q^6XP=\varepsilon(6hq^6XP+sq^{n_{11}}L^3_5+s^2(-q^{n_{21}}L^3_4PX-
q^{n_{22}}L^2_2X^3-\\
q^{n_{23}}L_2^2X^2P^{-1})+s^3
(3q^{n_{31}}L_1^1PX-q^{n_{32}}L^1_1P+q^{n_{33}}L^3_3PX)+
s^4(q^{n_{41}}L_1^2P+\\q^{n_{42}}X^{-2}L_1^2PX^{-1}-q^{n_{43}}L^3_2XP)+s^5
q^{n_{51}}L_1^3XP),
\end{multline*}
for the case $\ell=6$. In the last two formulas $m_{ij},n_{ij}$
are integers, whose exact values play no
role in the proof.
\end{proof}

\begin{proof}[Proof of Theorem~\ref{univ}] The second Hochschild
cohomology $H^2(\mathbb{Z}_\ell\ltimes D_q)=\CC^{r+1}$ is the
tangent space to the moduli space of all deformations. The
deformations coming from the family
$\{H({t},q)\}_{{t}\in (\CC^*)^r,q\in\CC^*}$ yield a
subspace in this tangent space. Lemma \ref{nontriv}
shows that this subspace is of dimension $r+1$, i.e., coincides
with the entire tangent space. This implies that the above family
furnishes a universal deformation. 
\end{proof}

\section{Generalized double affine Hecke algebras over $\Bbb C$}

\subsection{Definition}
Now we let $t_{kj}$ be complex numbers such that $\prod_j t_{kj}=1$
for all $k$. Such collections of numbers form
an algebraic torus $\Bbb T$. For $t\in \Bbb T$, define
an algebra $H(t,q)$ in the same way as $\widehat {\bold H}(q)$,
except that this is now an algebra over $\Bbb C$.
That is, the algebra $H(t,q)$ is generated over $\Bbb C$ by
$T_i,i=1,...,m$, with defining relations (\ref{defii}).
This family of algebras can be obtained by specializations
of a single (obviously defined) algebra ${\bold H}(q)$ over $\Bbb C[\Bbb T]$, in which $t_{kj}$
are central elements, and, yet more universally, of the algebra 
$\bold H$ over ${\mathcal R}:=\Bbb C[\Bbb T]\otimes \Bbb C[q,q^{-1}]$, in
which both $q$ and $t_{kj}$ are central elements. 

\subsection{The length filtration and PBW theorem}
We will now introduce an important length filtration on 
the algebras $\bold H$, $\bold H(q)$, and $H(t,q)$. To do so, let us  
note that the group $G$ is the group $W_+$ of even elements 
of an affine Weyl group $W$ of rank $m$, of types $\hat A_1\times \hat A_1$,
$\hat A_2$, $\hat B_2$, and $\hat G_2$, if $G$ is of types $D_4,
E_6, E_7, E_8$, respectively.\footnote{Here, by rank we mean the
number of nodes of the affine Dynkin diagram}

The group $W$ is generated by $s_1,...,s_{m}$ with the defining
relations $s_i^2=1$ and $(s_ks_{k+1})^{d_k}=1, k\in \Bbb Z_m$.
The isomorphism of $\eta_0: G\to W_+$ is given by the formula $T_k\to
s_ks_{k+1}$. 

Let $\sigma: G\to G$ be the automorphism defined by the formulas 
$\sigma(T_1)=T_1^{-1}$, $\sigma(T_2)=T_2^{-1}$, 
and in the $D_4$ case $\sigma(T_3)=T_2T_3^{-1}T_2^{-1}$. 
It is easy to see that the homomorphism 
$\eta_0$ extends to an isomorphism 
$\eta: \lbrace{1,\sigma\rbrace}\ltimes G\to W$, 
given by the formula $\eta_0(\sigma)=s_2$. 

Let us now construct a deformation of $\eta_0$. Let 
$f\mapsto \overline{f}$ be the automorphism of ${\mathcal R}$ defined by the formula 
$\overline{q}=q^{-1}, \overline{t_{kj}}=t_{k,-j}^{-1}$
(where $-j$ is taken modulo $d_k$).
Let the algebra $\bold A$ be generated by ${\mathcal R}$ and 
additional generators $s_k$ with defining relations 
$$
s_k^2=1,\quad s_kf=\overline{f}s_k,\ f\in {\mathcal R}, 
$$
and 
$$
\prod_{j=1}^{d_k}(s_ks_{k+1}-q^{-\delta_{km}}t_{kj})=0.
$$
It is clear that we have a homomorphism of algebras 
$\eta: \bold H\to \bold A$ 
defined by the formula 
$\eta(T_k)=q^{\delta_{km}}s_ks_{k+1}$. 

Let $\sigma: \bold H\to \bold H$ be the automorphism of algebras defined 
by the formulas $\sigma(f)=\overline{f}$ for $f\in {\mathcal R}$, 
$\sigma(T_1)=T_1^{-1}$, $\sigma(T_2)=T_2^{-1}$, 
and in the $D_4$ case $\sigma(T_3)=T_2T_3^{-1}T_2^{-1}$
(clearly, $\sigma^2=1$). It is easy to see that the homomorphism 
$\eta$ extends to an isomorphism 
$\eta: (\Bbb C[\sigma]/(\sigma^2=1))\ltimes \bold H\to \bold A$, 
given by the formula $\eta(\sigma)=s_2$. 
To see that this is an isomorphism, it suffices to construct the
inverse, which is given by the formulas 
$\eta^{-1}(s_2)=\sigma$, $\eta^{-1}(s_1)=T_1\sigma$, 
$\eta^{-1}(s_3)=\sigma T_2$, and for type $D_4$, $\eta^{-1}(s_4)=\sigma
T_2T_3$. 

Define the filtration $F_L^\bullet$ on $\bold A$ and $\bold H$ 
by the formulas $\deg({\mathcal R})=0$, $\deg(s_i)=1$. We call this
filtration {\it the length filtration}. 
 
For any element $x\in W$, fix a reduced decomposition $w(x)$ of
$x$. Also, for any word $w$ in the letters $s_i$, 
define the element $T_w\in \bold A$ 
to be the product of the generators $s_i$ according to the
word $w$. 

The next theorem (which can also be found in \cite{ER}) 
is a PBW theorem for ${\bold H}$.
It is formulated in terms of the length filtration.

\begin{theorem}\label{pbwl}
(i) The elements $T_{w(x)}$ form a basis of $\bold A$ as a 
left ${\mathcal R}$-module. The elements 
$T_{w(x)}$ for even $x$ form a basis of $\bold H$ over ${\mathcal R}$.  

(ii) The elements $T_{w(x)}$ with ${\rm length}(x)\le N$
form a basis of $F_L^N\bold A$ as a 
left ${\mathcal R}$-module. The elements 
$T_{w(x)}$ for even $x$ of length $\le N$ 
form a basis of $F_L^N\bold H$ over ${\mathcal R}$.  

(iii) The ${\mathcal R}$-modules $F_L^N \bold H/F_L^{N-1}\bold H$, 
$F_L^N\bold A/F_L^{N-1}\bold A$ are free.
The Hilbert series of $\bold H,\bold A$ 
under the length filtration is the same as those of $W,W_+$. 
\end{theorem}

\begin{proof}
It is clear that (i) and (iii) follow from (ii), so it
suffices to prove (ii).

First of all, the elements $T_{w(x)}$ are linearly independent
after reduction modulo the ideal $t_{kj}=e^{2\pi {\rm i}j/d_k}$, 
so by Theorem \ref{pbw} they are linearly independent in 
$\bold A$. It remains to show that $T_{w(x)}$ 
with length of $x$ being $\le N$ is a spanning set of $F_L^N
\bold A$. 

Let us write the relation
$$
\prod_{j=1}^{d_k}(s_ks_{k+1}-q^{-\delta_{km}}t_{kj})=0
$$
as a deformed braid relation:
$$
s_ks_{k+1}...+S.L.T.=t_{k}s_{k+1}s_k...+S.L.T.,
$$
where $t_{k}:=(-1)^{d_k+1}q^{-\delta_{km}d_k}t_{k1}...t_{kd_k}$, S.L.T.
mean ``smaller length terms'',
and the products on both sides have length $d_k$.
This can be done by multiplying the relation by $s_ks_{k+1}...$
($d_k$ factors).

Clearly, $T_w$ for all words $w$ of length $\le N$ span $F_L^N \bold A$.
So we just need to take any word $w$ of length $\le N$ and
express $T_w$ via $T_{w(x)}$ for $x\in W$, ${\rm length}(x)\le
N$. 

It is well known from the theory of Coxeter groups (see
e.g. \cite{B}) that
using the braid relations, one can turn any non-reduced word
into a word that is not square free, and any
reduced expression of a given element
of $W$ into any other reduced expression of the same element.
Thus, if $w$ is non-reduced, then by using the deformed braid
relations and the relations $s_i^2=1$, 
we can reduce $T_w$ to a linear combination of $T_u$ with words
$u$ of smaller length than $w$.
On the other hand, if $w$ is a reduced expression for some
element $x\in W$, then using the deformed braid relations
we can reduce $T_w$ to
a linear combination of $T_u$ with $u$ shorter than $w$, and
$T_{w(x)}$.
Thus $T_{w(x)}$ are a spanning set.
The theorem is proved.
\end{proof}

\begin{corollary} The Gelfand-Kirillov dimension of ${\bold H}(q)$
as an algebra over $\Bbb C[\Bbb T]$, and of $H(t,q)$ is 2.
\end{corollary}

\subsection{The general properties of $H(t,q)$ when $q$ is a root
of unity}

Recall that an algebra $A$ is PI of degree $K$ if all
polynomial identities of the matrix algebra of size $K$ are
satisfied in $A$.

\begin{theorem}\label{pol}
Let $q$ be a root of unity such that $q^\ell$ has degree $N$.
Then $H(t,q)$ and ${\bold H}(q)$ are PI algebras of degree
$\le \ell N$.
\end{theorem}

\begin{proof}  According to the proof of
Theorem \ref{pbw}, the algebra $B(q)$ has an
embedding into ${\rm Mat}_{\ell N}(R)$,
where $R$ is the ring of regular functions
on a formal 2-dimensional polydisk
(obtained by considering the formal neighborhood
of an irreducible $\ell N$-dimensional representation of $B(q)$).
It is shown in the proof of Theorem \ref{pbw} that this embedding can
be deformed into an embedding of $\widehat {\bold H}(q)$ into
${\rm Mat}_{\ell N}(R)[[\tau]]$. This implies the result
for the formal algebra $\widehat {\bold H}(q)$.

Now let us establish the result for the algebra ${\bold H}(q)$.
Theorems \ref{pbw},\ref{pbwl} imply that the algebra
$\widehat {\bold H}(q)$ is the formal completion
of ${\bold H}(q)$ with respect to the ideal
defined by the equation $t=1$, and the
natural map ${\bold H}(q)\to \widehat {\bold H}(q)$ is injective.
This implies the result for ${\bold H}(q)$ and hence for $H(t,q)$ for
special $t$.
\end{proof}

Let $U=U_N\subset \Bbb T$ be the Zariski open set,
defined by the following condition: one cannot choose
nonnegative integers $p_{kj}$ with 
$\sum_j p_{kj}=p< \ell N$ for all $k=1,...,m$
such that $\prod_{k,j}(u_{kj})^{p_{kj}}=q^p$.
It is clear that $U$ is nonempty.

\begin{theorem}\label{dim}
Let $q$ be a root of unity such that $q^\ell$ has degree
$N$. Then every irreducible representation
of $H(t,q)$ has dimension
$\le\ell N$. For $t\in U$, the dimension is exactly $\ell N$.
\end{theorem}

\begin{proof} Let $V$ be a nonzero finite dimensional
representation of $H(t,q)$ of dimension $p<\ell N$.
Then the product of determinants of $T_k$ is equal to $q^p$.
Thus $q^p$ should equal the product of $u_{kj}$ with
multiplicities. For $t\in U$, this leads to a contradiction, so
for such $t$ we have $p\ge \ell N$.

It remains to prove the opposite inequality
for irreducible representations $V$.
This follows from Theorem \ref{pol} and
the following well known theorem
is due to Kaplansky.

\begin{theorem}
If $A$ is a PI algebra of degree $\le K$ then
any irreducible $A$-module is finite dimensional
and has dimension $\le K$.
\end{theorem}
\end{proof}

{\bf Remark.} Note that if $q$ is not a root of unity, then
$H(t,q)$ does not have finite dimensional representations for
generic $t$. This can be seen by taking the determinant of the
relation $\prod_k T_k=q$.

\subsection{The algebras $H(t,q)$
when $q$ is a root of unity, for generic $t$}

\begin{theorem}\label{Azu}
Let $q$ be a root of unity such that $q^\ell$ has degree $N$.
Then for $t\in U$,
the algebra $H(t,q)$ is an Azumaya algebra of rank $\ell N$
over an affine 2-dimensional $\Bbb C$-scheme $S(t,q)$ of finite type.
\end{theorem}

\begin{proof}
We first recall the following theorem
about PI algebras.

\begin{theorem} (M. Artin, see \cite{Ja})
Let $A$ be a finitely generated algebra over $\Bbb C$, which is
PI of degree $\le K$. If the dimension of every irreducible 
finite dimensional $A$-module
is $K$ then $A$ is an Azumaya algebra,
and its center is finitely generated.
\end{theorem}

This result and Theorem \ref{dim} implies that
for $t\in U$, the algebra $H(t,q)$ is Azumaya, with finitely
generated center $Z(t,q)$.

Let $S(t,q)$ be the spectrum of $Z(t,q)$.
Since the Gelfand-Kirillov dimension of $H(t,q)$ is 2,
$\dim S(t,q)=2$.

\end{proof}

\begin{corollary} Let $q=1$, $t\in U=U_1$, and the minimal polynomial
of the generator $c=T_3$ has a simple root. Then $H(t,q)$  
is the endomorphism algebra of a vector bundle over $S(t,q)$.
\end{corollary}

\begin{proof}
Let $e$ be the projection to the eigenspace of $c$
corresponding to the simple root of the minimal polynomial 
(it is a polynomial of degree $\ell$ of the generator $c=T_3$). 
If $q=1$ and $t\in U$
then $H(t,q)$  is an Azumaya algebra, so  
the map $Z(t,q)\to eH(t,q)e$ given by $z\to ze$ is an
isomorphism, and $H(t,q)e$ is a projective module over $eH(t,q)e$. This gives
rise to the required vector bundle.
\end{proof}

\subsection{Smoothness of $S(t,q)$}

\begin{theorem}\label{smooth} If $q$ is a root of unity and $t\in
U$ with $u_{kj}\ne u_{kj'}$ for $j\ne j'$
then $S(t,q)$ is a smooth surface.
\end{theorem}

\begin{proof} Let $K=\ell N$.
Since $H(t,q)$ is an Azumaya algebra,
$S(t,q)$ is isomorphic to the moduli space $M(t,q)$ of irreducible
($K$-dimensional) representations of $H(t,q)$. The space $M(t,q)$
is a quotient of the affine variety of $K$-dimensional
matrix representations of $H(t,q)$ by the free $PGL_K$-action, so
it is an affine variety.

More specifically, each connected component\footnote{A priori, the
connected components of $M(t,q)$ correspond to various
multiplicities of eigenvalues of $T_j$. However, we will see
later in Proposition \ref{irraf} 
that in fact the varieties $M(t,q)$ and $S(t,q)$ are connected,
i.e. eigenvalues of $T_j$ always occur with equal multiplicities.} of the 
the space $M(t,q)$ is a quotient
by the free $PGL_K$-action of the
subvariety $Y$ in the product of $m$ special semisimple conjugacy classes
$C_1,...,C_m$ in $GL_K$ (those of $T_j$)
defined by the equation $T_1...T_m=q$. Thus we
just need to show that this subvariety is smooth.

Let us first deal with the cases $E_6,E_7,E_8$. In these cases,
let us fix the matrix $c$ to be diagonal and consider the map
$\mu: C_1\times C_2\to GL_K$, $\mu(a,b)=ab$. It lands
in a fixed coset of $SL_K$ in $GL_K$. Our job is to show that
$qc^{-1}$ is a regular value of $\mu$. To do so,
it suffices to show that the differential of $\mu$:
$d\mu(a,b)=(da)b+a(db)$, is surjective onto $sl_K$.
Let $v=[x,a],w=[y,b]$ be two tangent vectors to $C_1,C_2$
at $a,b$, respectively.
Then we get
$$
d\mu(a,b)(v,w)=[x,a]b+a[y,b]
$$
To show that the map $x,y\to [x,a]b+a[y,b]$ is surjective, let
us assume that $z\in gl_K$ is orthogonal to its image.
Then we have ${\rm Tr}([x,a]bz)=0,{\rm Tr}(a[y,b]z)=0$ for all
matrices $x,y$. So $[a,bz]=[b,za]=0$. Set $u=bza$. Then
$[a,u]=[b,u]=0$. Hence $u$ is a scalar (since $a,b$ is an
irreducible collection of matrices), and $z=\lambda
b^{-1}a^{-1}$, where $\lambda$ is a scalar.
So the map in question is surjective (onto the tangent space
of the $SL_K$-coset).

The case of $D_4$ is similar. In this case we have a map
$\nu: C_1\times C_2\times C_4\to GL_K$ given by
$\nu(a,b,d)=dab$. We need to show that the value $qc^{-1}$ is
regular. The computation of the orthogonal complement of the
image of the differential of $\nu$ will give equations
$[d,abz]=0$, $[a,bzd]=0$, $[b,zda]=0$.
Set $u=bzda$. Then $[a,u]=[b,u]=[d,u]=0$. Since $a,b,d$ are
irreducible, $u$ is a scalar, and $z=\lambda
b^{-1}a^{-1}d^{-1}$,
as desired.
\end{proof}

\section{The geometry of the algebras $H(t,q)$}

\subsection{Other filtrations}

In this subsection we will study filtrations 
on $\bold H(q)$ which are, in reality, simpler and more useful 
that $F_L$, but whose study requires computer calculations.
They will be used to study the geometry of the algebras 
$H(t,q)$, in particular their intimate connection with Del Pezzo
surfaces. To be more specific, these filtrations correspond to 
the simplest compactifications of affine del Pezzo surfaces,
which are described below. 

Namely, we introduce increasing filtrations
$F_{jkl}^\bullet {\bold H}(q)$ on ${\bold H}(q)$ as follows:
$\deg \Bbb C[\Bbb T]=0$, $\deg(c)=0$, and

$F_{111}^\bullet$ for $D_4$: $\deg(b)=1$, $\deg(d)=1$, $\deg(a)=1$;

$F_{112}^\bullet$ for $D_4$: $\deg(b)=1$, $\deg(d)=1$;

$F_{111}^\bullet$ for $E_6$: $\deg(b)=1$, $\deg(b^2)=1$, $\deg(b^2cb)=1$;

$F_{112}^\bullet$ for $E_6$: $\deg(b)=1$, $\deg(b^2)=1$;

$F_{123}^\bullet$ for $E_6$: $\deg(b)=1$;

$F_{112}^\bullet$ for $E_7$: $\deg(b)=1$, $\deg(b^2)=1$;

$F_{123}^\bullet$ for $E_7$: $\deg(b)=1$;

$F_{123}^\bullet$ for $E_8$: $\deg(b)=1$.

Since these elements are generators, these filtrations are well defined.
Similarly one defines a filtration $F_{jkl}^\bullet H(t,q)$ for specific
points $t$.

{\bf Remark.} The subscripts are the degrees 
of the generators of the center of $H(t,1)$ with respect to the
corresponding filtration; this will be discussed below. 

\begin{theorem}\label{pbw1} (The PBW theorem)
The spaces $F_{jkl}^N{\bold H}(q)/F_{jkl}^{N-1}{\bold H}(q)$ 
are free $\Bbb C[\Bbb T]$-modules.
\end{theorem}

\begin{proof}
 Theorem \ref{pbw1} is proved in 
subsections \ref{pbws1},\ref{pbws2}.
The structure of the proof is as follows.

Let $F=F_{jkl}$. Let $m_r$ be the dimension of $F^r H(1,q)$, and 
$\lbrace{g_s, s\ge 1\rbrace}$ be a labeling of elements of $G$
by positive integers, such that $g_1,...,g_{m_r}$ is a basis of
$F^rH(1,q)$. 

We will find a collection of ``legal'' monomials 
$\lbrace{h_s\in {\bold H}(q), s\ge 1\rbrace}$ in generators,
which satisfy the following conditions: 

(i) If $m_{r-1}<s\le m_r$ then $h_s\in F^r{\bold H}(q)$, and 
$h_s$ with $s\le m_r$ span $F^r{\bold H}(q)$ 
over $\Bbb C[\Bbb T]$;

(ii) $h_s$ specialize to $g_s$ at $t=1$.  

Property (i) implies that 
for $m_{r-1}<s\le
m_r$ the images $h_s'$ of $h_s$ \linebreak  
in $F^r{\bold H}(q)/F^{r-1}{\bold H}(q)$ span this module
over $\Bbb C[\Bbb T]$. Thus, it remains to show that 
they are linearly independent. 

To do so, we assume that we have a nontrivial linear relation
$$
\sum_{s=m_{r-1}+1}^{m_r}f_s(t)h_s'=0,\ f_s\in \Bbb C[\Bbb T].
$$ 
Then we have a linear relation 
$$
\sum_{s=1}^{m_r}f_s(t)h_s=0.
$$
Expanding this relation in a power series near $t=1$, we get 
a linear relation in $\widehat{\bold H}(q)$. 
This relation is nontrivial because of condition (ii). 
Thus we obtain a contradiction with Theorem \ref{pbw}.
This completes the proof. 
\end{proof} 

In the case $t=1$, the filtrations $F_{jkl}^\bullet$ are very easy to
understand. In this case, we have a natural (up to scaling by
powers of $q$) basis of $H(t,q)=B(q)$ corresponding to group
elements, and it is easy to see that our filtrations
are compatible with this basis. So let us say which
basis elements have degree $\le N$.

Realizing $G$ as $\Bbb Z_\ell\ltimes \Bbb Z^2$,
we can write any group element as a product $c^jY$, where $0\le
j\le \ell-1$ and $Y$ belongs to the lattice.
The basis elements of degree $\le N$ are then those products
$c^jY$ for which $||Y||\le N$, where $||Y||$ is a certain norm on
$\Bbb R^2$, depending on a particular filtration.
We will describe these norms in all cases.

$D_4$, $F_{111}$. The lattice is hexagonal (generated by two vectors
$v,w$ of equal length making angle $60^o$ with each other)
and the norm is such that the unit ball is the
hexagon whose vertices are $v$ and its images under $\Bbb
Z_6$.

$D_4$, $F_{112}$.
The lattice is rectangular, and 
the norm is $|x|+|y|$.

$E_6$. The lattice is hexagonal, and

(i) For $F_{111}$, the norm is such that the unit ball
is the triangle whose vertices are
$v+w$ and its two images under $\Bbb Z_3$.

(ii) For $F_{112}$, the norm is such that the unit ball is the
hexagon whose vertices are $v$ and its images under $\Bbb
Z_6$.

(iii) For $F_{123}$, the norm is such that the unit ball is the
triangle whose vertices are $v$ and its images under $\Bbb
Z_3$.

$E_7$. The lattice is rectangular, and 

(i) For $F_{112}$, the norm is ${\rm
max}(|x|,|y|)$.

(ii) For $F_{123}$, the norm is $|x|+|y|$.

$E_8$, $F_{123}$. The lattice is hexagonal, 
and the norm is such that the unit ball is the
hexagon whose vertices are $v$ and its images under $\Bbb
Z_6$.

This implies the following result.

\begin{proposition}\label{hilser}
The Poincar\'e series
of $H(t,q)$ and ${\bold H}(q)$ with respect to the filtration
$F^\bullet_{jkl}$ is:
$$
P_{D_4,111}(z)=2(1+\frac{6z}{(1-z)^2}),
$$
$$
P_{D_4,112}(z)=2(1+\frac{4z}{(1-z)^2}),
$$
$$
P_{E_6,111}(z)=3(1+\frac{9z}{(1-z)^2}),
$$
$$
P_{E_6,112}(z)=3(1+\frac{6z}{(1-z)^2}),
$$
$$
P_{E_6,123}(z)=3(1+\frac{3z}{(1-z)^2}),
$$
$$
P_{E_7,112}(z)=4(1+\frac{8z}{(1-z)^2}),
$$
$$
P_{E_7,123}(z)=4(1+\frac{4z}{(1-z)^2}),
$$
$$
P_{E_8,123}(z)=6(1+\frac{6z}{(1-z)^2}),
$$
\end{proposition}

\subsection{The associated graded algebras of ${\bold H}(q)$ and
$H(t,q)$}

Let ${\bold H}_0^{jkl}(q)$ and $H_0^{jkl}(t,q)$ be the 
associated graded algebras
of ${\bold H}(q)$ and $H(t,q)$ with respect
to the filtration $F_{jkl}$.

Let us give a description of some of the algebras $H_0^{jkl}(t,q)$ by
generators and relations. A similar description (with $t$ being
variables) is valid for ${\bold H}_0^{jkl}(q)$. 

For $\ell\ge 3$ and a monic polynomial $p$ of degree $\ell$, define 
an algebra $K_\ell(p)$ to be the free algebra generated by  
elements $c,d$ modulo the ideal generated by elements 
\begin{equation}\label{graded}
p(c),\
d c^2 d,\ d c^3 d,\dots,\ d c^{\ell-2} d,\
d c d - q^{-\ell} c d c^{\ell-1} d c,\
d c d c d.
\end{equation}

\begin{theorem}\label{genrel}
In the cases $E_6,E_7,E_8$, the algebra ${\bold H}_0^{123}(t,q)$ 
is isomorphic to the algebra $K_\ell(p)$, 
where $p(x)=
\prod_{j=1}^\ell(x-u_{3j})$, and $\ell=3$, $4$, $6$ for
$E_6$, $E_7$, $E_8$ respectively.
\end{theorem}

\begin{proof}
The proof is given in subsection \ref{assgr}.
Note that in subsection \ref{assgr} we obtain 
the same relations with $-q^{-\ell}$ replaced with 
$(-1)^\ell q^{-\ell}$, due to the fact that generators have been
rescaled. 
\end{proof}

A similar description of the graded algebra exists in the $D_4$
case. Namely, we have 

\begin{proposition}\label{genreld4}
In the case $D_4$, the algebra $H_0^{112}(t,q)$ 
is generated by $c,z_1,z_2$ with defining relations
$$
p(c)=0,\ z_1cz_1=0,\ z_2cz_2=0,\
cz_1cz_2c=q^{2}z_2cz_1,\ z_1z_2=q^{-2}z_2c^2z_1.
$$ 
On the other hand, the algebra $H_0^{111}(t,q)$ for $D_4$ 
is generated by $c,w_1,w_2,w_3$ with defining relations
$$
p(c)=0,\
w_1 w_3=0,\ w_1 w_2=0\ ,w_2 w_3=0,
$$
$$
w_3 c^2 w_1=0,\ w_2 c^2 w_1=0,\ w_3 c^2 w_2=0,
$$
$$
w_1 c w_1=0,\ w_2 c w_2=0,\ w_3 c w_3=0,
$$
$$
w_3 c w_1 = q^2 c w_1 c w_3 c,\
w_2 c w_1 = q^{-2} c w_1 c w_2 c,\
w_3 c w_2 = q^{-2} c w_2 c w_3 c.
$$
\end{proposition}

\begin{proof}
The proof is analogous to the proof of Theorem \ref{genrel},
using the presentations of $H(t,q)$ in Section 8; 
the elements $z_1$ and $z_2$ are the images in the graded algebra of
the elements $c^{-1}b$ and $dc^{-1}$, while 
the elements $w_1,w_2,w_3$ are the images 
of $c^{-1}d^{-1},db,b^{-1}c^{-1}$.  
\end{proof}

These results show that the algebra $H_0^{jkl}(t,q)$ 
for the considered filtrations does not depend on
$t_{ij}$ with $i\ne 3$. 

{\bf Remark.} These results imply that the PI degree 
in Theorem \ref{pol} is exactly $\ell N$ (since it is 
so for the corresponding associated graded algebras).

\subsection{The geometric characterization of the 
associated graded algebras.} 

To characterize the associated graded algebras geometrically, 
let us recall the theory of
noncommutative curves (\cite{SV}).

Let $X$ be a projective algebraic curve,
$\sigma$ an automorphism of $X$, and ${\mathcal L}$ an ample
line bundle on $X$. Then one can define the twisted
homogeneous coordinate ring $B(X,\sigma,{\mathcal L})$ as follows.
This is a $\Bbb Z_+$-graded ring, and
$B(X,\sigma,{\mathcal L})[n]={\rm Hom}({\mathcal O},
\otimes_{j=0}^{n-1}{\mathcal L}^{\sigma^j})$. The multiplication is defined
by the formula $a*b=\sigma_*^{{\rm deg} b}(a)\otimes b$
for homogeneous $a,b$. In noncommutative algebraic geometry,
this twisted homogeneous coordinate ring
is viewed as the homogeneous coordinate ring of a noncommutative projective
curve.

Similarly, if $E$ is a vector bundle on $X$
equivariant under $\sigma$ then one can define
the graded algebra $B(X,\sigma,{\mathcal L},E)$ by
$$
B(X,\sigma,{\mathcal L},E)[n]={\rm Hom}(E,
E\otimes (\otimes_{j=0}^{n-1}{\mathcal L}^{\sigma^j})).
$$
with multiplication as above.

Let $X_n'=\Bbb P^1_1\cup \Bbb P^1_2\cup...\cup \Bbb P^1_n$ 
be a chain of projective lines, i.e. the point 
$0$ of $\Bbb P^1_i$ is identified with the point $\infty$ 
of $\Bbb P^1_{i+1}$ for $1\le i\le n-1$. 

Let $X_n$ be the union of $n$ projective
lines forming an $n$-gon (i.e. each two
consecutive ones in a cyclic order intersect at a point).
Clearly, $X_n$ is
obtained from $X_n'$ by gluing the point $\infty$ of 
$\Bbb P^1_1$ with the point $0$ of $\Bbb P^1_n$. 
 For $n=1$, $X_n$ is a single $\Bbb P^1$ with a node, and we will denote
it simply by $X$.

Let $P_i$ be a smooth point lying on the $i$-th component of $X_n$.
The group $\Bbb C^*$ acts naturally on $X_n$;
let $\sigma$ be the action of $q^\ell\in \Bbb C^*$.
Set ${\mathcal L}={\mathcal O}(\sum_i P_i)$,

Let $p_t(x)=x^\ell+\alpha_1x^{\ell-1}+...+\alpha_\ell$ be the 
monic polynomial of degree $\ell$ annihilating $c$
($p_t(x)=\prod_j(x-u_{3j})$). Let $A_t$ be the companion matrix 
corresponding to the polynomial $p_t$. This is the
$\ell$-by-$\ell$ matrix defined by the formula 
$$
A_tv_i=v_{i+1}, 1\le i<\ell,\
A_tv_\ell=-\alpha_1v_{\ell}-\alpha_2v_{\ell-1}-...-\alpha_\ell v_1,
$$  
where $v_i$ is the standard basis of $\Bbb C^\ell$;
thus $p_t(A_t)=0$. 

Consider the trivial vector bundle of rank $\ell$ on
$X_n'$. Let $V_0$ be its fiber at $0\in \Bbb P^1_n$, and 
$V_\infty$ its fiber at $\infty\in \Bbb P^1_1$. 

Let $A$ be an invertible $\ell$ by $\ell$-matrix. 
Let $E(A)$ be the vector bundle on $X_n$
obtained from the trivial bundle of rank $\ell$ on 
$X_n'$ by gluing the fibers $V_0,V_\infty$ using the map
$A: V_\infty\to V_0$. Thus, if $\Delta$ is an effective divisor 
on $X_n$ not containing the gluing points then 
sections of $E(A)$ with poles at $\Delta$ are 
collections of $\Bbb C^\ell$-valued
rational functions $\phi_1,...,\phi_n$ of one variable
$z$ with poles at $\Delta$ which satisfy the conditions
$$
\phi_1(0)=\phi_2(\infty),...,\phi_{n-1}(0)=\phi_n(\infty),\
\phi_n(0)=A\phi_1(\infty).
$$
Obviously, the bundle $E(A)$ is $\Bbb C^*$-equivariant.

\begin{theorem} \label{assgrad}
(i) The algebra $H_0^{jkl}(t,q)$ is isomorphic
to $B(X_n,\sigma,{\mathcal L},E(A_t^\varepsilon))$,
where: 

for $D_4$ and $jkl=111$, $n=3$, $\varepsilon=-1$; 

for $D_4$ and $jkl=112$, $n=2$, $\varepsilon=-1$; 

for $E_6,E_7,E_8$, $jkl=123$, $n=1$, $\varepsilon=1$. 

(ii) Let $u_{31}$ be a simple root of the minimal polynomial 
of $c$. Let $e$ be the projector to the
$u_{31}$-eigenspace of $c$ in $\Bbb C[c]\subset H_0^{jkl}(t,q)$.
Then the ``spherical subalgebra'' $eH_0^{jkl}(t,q)e$
is isomorphic to $B(X_n,\sigma,{\mathcal L})$ (for $n$ as above).
\end{theorem}

\begin{proof} The second statement follows from the first one, so
it suffices to prove (i). Let us do it first in the cases 
$E_6,E_7,E_8$. 

We start with an explicit description of the algebra 
$B(X,\sigma,{\mathcal L},E(A_t))$. 

By the definition, 
$B(X,\sigma,{\mathcal L},E(A_t))[n]$ is the space of 
rational functions $f$ of one variable $z$ with values in ${\rm
Mat}_\ell(\Bbb C)$ with divisor of poles dominated by the divisor 
$(1)+(q^\ell)+...+(q^{(n-1)\ell})$, and such that 
$f(0)=A_tf(\infty)A_t^{-1}$. Furthermore, the multiplication law 
in $B(X,\sigma,{\mathcal L},E(A_t))$ is given by the formula
$(f*g)(z)=f(q^{-\ell m}z)g(z)$ for homogeneous elements $f,g$
such that $g$ has degree $m$. Thus, for example, 
$B(X,\sigma,{\mathcal L},E(A_t))[0]$ has basis
$1,A_t,A_t^2,...,A_t^{\ell-1}$, while 
$B(X,\sigma,{\mathcal L},E(A_t))[1]$ is the space of functions of
the form
$$
f(z)=\frac{Mz-A_tMA_t^{-1}}{z-1},
$$
where $M$ is any matrix of size $\ell$ by $\ell$. 

Now we will define a homomorphism 
of graded algebras $\xi: H_0^{123}(t,q)\to B(X,\sigma,{\mathcal L},E(A_t))$. 
It is defined by the formulas
$$
\xi(c)=A_t,\ \xi(d)=\frac{E_{11}z-A_tE_{11}A_t^{-1}}{z-1},
$$
where $E_{ij}$ is an elementary matrix. 
It is straightforward to check that 
the relations of $H_0^{123}(t,q)$ given in Theorem \ref{genrel} 
are satisfied, hence $\xi$ is well defined. 
Moreover, it is easy to check that 
$B(X,\sigma,{\mathcal L},E(A_t))$ is generated by degrees $0$ and
$1$, which implies that $\xi$ is surjective. By comparing the
Poincar\'e series we find that $\xi$ is bijective. 
This proves (i) for $E_6,E_7,E_8$. 

Now let us handle the case $D_4$, the 112 filtration.
By the definition, 
$B(X_2,\sigma,{\mathcal L},E(A_t^{-1}))[n]$ is the space of pairs
rational functions $(f_1,f_2)$ of one variable $z$ with values in ${\rm
Mat}_2(\Bbb C)$ with divisor of poles dominated by the divisor 
$(1)+(q^2)+...+(q^{2(n-1)})$, and such that $f_1(0)=f_2(\infty)$,
$A_tf_2(0)A_t^{-1}
=f_1(\infty)$. Thus 
$B(X_2,\sigma,{\mathcal L},E(A_t^{-1}))[0]$ has basis
$(1,1),(A_t,A_t)$, while 
$B(X_2,\sigma,{\mathcal L},E(A_t^{-1}))[1]$ is the space of pairs of
functions 
$(f_1,f_2)$ of
the form
$$
f_1(z)=\frac{A_tM_1A_t^{-1}z-M_2}{z-1},\
f_2(z)=\frac{M_2z-M_1}{z-1},
$$
where $M_i$ are any $2$ by $2$ matrices. 

Now we will define a homomorphism 
of graded algebras $\xi: H_0^{112}(t,q)\to B(X_2,\sigma,{\mathcal
L},E(A_t^{-1}))$. 
It is defined by the formulas
$$
\xi(c)=A_t,\ \xi(z_1)=(-\frac{E_{11}}{z-1},
\frac{E_{11}z}{z-1}),\
\xi(z_2)=(\frac{A_tE_{11}A_t^{-1}z}{z-1},
\frac{-E_{11}}{z-1}).
$$
Similarly to the $E_{6,7,8}$ cases, $\xi$ is well defined 
and is an isomorphism. 

Finally, we consider the 111 filtration for $D_4$. 
By the definition, \linebreak
$B(X_3,\sigma,{\mathcal L},E(A_t^{-1}))[n]$ is the space of triples of
rational functions $(f_1,f_2,f_3)$ of one variable $z$ with values in ${\rm
Mat}_2(\Bbb C)$ with divisor of poles dominated by the divisor 
$(1)+(q^2)+...+(q^{2(n-1)})$, and such that $f_1(0)=f_2(\infty)$,
$f_2(0)=f_3(\infty)$, $A_tf_3(0)A_t^{-1}=f_1(\infty)$. 
Thus  
$B(X_3,\sigma,{\mathcal L},E(A_t^{-1}))[0]$ has basis
$(1,1,1),(A_t,A_t,A_t)$, while 
$B(X_3,\sigma,{\mathcal L},E(A_t^{-1}))[1]$ is the space of triples of
functions 
$(f_1,f_2,f_3)$ of
the form
$$
f_1(z)=\frac{A_tM_1A_t^{-1}z-M_3}{z-1},
f_2(z)=\frac{M_2z-M_1}{z-1},
f_3(z)=\frac{M_3z-M_2}{z-1},
$$
where $M_i$ are any $2$ by $2$ matrices. 

Now we will define a homomorphism 
of graded algebras $\xi: H_0^{111}(t,q)\to B(X_3,\sigma,{\mathcal
L},E(A_t^{-1}))$. 
It is defined by the formulas
$$
\xi(c)=A_t, \xi(w_1)=(-\frac{E_{11}}{z-1},
\frac{E_{11}z}{z-1},0),
$$
$$
\xi(w_2)=(\frac{A_tE_{11}A_t^{-1}z}{z-1},0,
\frac{-E_{11}}{z-1}),\ \xi(w_3)=(0,-\frac{E_{11}}{z-1},
\frac{E_{11}z}{z-1}).
$$
Similarly to the $E_{6,7,8}$ cases, $\xi$ is well defined 
and is an isomorphism. 
\end{proof}

\subsection{The center of $H(t,q)$}
We return to the study of the center $Z(t,q)$ and the scheme 
$S(t,q)$ when $q$ is a root of unity, such that $q^\ell$ has
order $N$. 

\begin{proposition}\label{irraf} For any $t$, 
the scheme $S(t,q)$ is an irreducible affine algebraic surface.
\end{proposition}

\begin{proof}
In the $D_4$ case, the theorem follows from \cite{Ob}.
Thus consider the cases $E_6,E_7,E_8$.
The associated graded algebra ${\rm gr}Z(t,q)$
is a subalgebra in the center $Z({\rm gr}H(t,q))$
of ${\rm gr}H(t,q)$.
As follows from our description of ${\rm gr}H(t,q)=H_0^{123}(t,q)$
in Theorem \ref{assgrad}, the center  $Z({\rm gr}H(t,q))$
is the function algebra of the cone of a nodal $\Bbb P^1$.
So this algebra has no zero divisors, which implies the result.
\end{proof}

\begin{theorem} \label{grcenter} Let $q$ be a root of unity
such that $q^\ell$ has order $N$. 
Consider the filtration $F_{jkl}$ on $H(t,q)$, where $jkl=111$
for $D_4$ and $jkl=123$ for $E_{6,7,8}$. 
Then for any $t$:

(i) ${\rm gr}Z(t,q)=Z({\rm gr}H(t,q))$ (where $Z(A)$ denotes
the center of an algebra $A$, and ${\rm gr}$ is taken with
respect to $F_{jkl}$);

(ii) The Poincar\'e series of $Z({\rm gr}H(t,q))$ is:
$$
Q_{D_4}(z)=1+\frac{3z^N}{(1-z^N)^2};
$$
$$
Q_{E_{6,7,8}}(z)=1+\frac{z^N}{(1-z^N)^2}.
$$
\end{theorem}

\begin{proof} Statement (ii) follows from the description of
${\rm gr}H(t,q)$given in
Theorem \ref{assgrad}.
So let us prove (i).

A priori, ${\rm gr}Z(t,q)$ is a subalgebra
in $Z({\rm gr}H(t,q))$. We must show that in fact, these two algebras
coincide.

The coincidence of the two algebras is easy in the case $t=1$, by a direct
computation. Therefore, it suffices to establish the coincidence of the
two algebras for generic (or even Weil generic) $t$.

In the $D_4$ case, it is easy to produce
three generators in degree $N$ of
${\rm gr}Z(t,q)$: one should
consider the Demazure-Lusztig realization of
$H(t,q)$ as reflection-difference operators on functions of $x$
(see \cite{Sa,St,Ob} for a definition) and take the element $Z_1=x^N+x^{-N}$
and its images under the $SL_2(\Bbb Z)$-action.
This implies the desired result for $D_4$, so we can now focus on
$E_6,E_7,E_8$.

Recall that for a $\Bbb Z_+$-filtered algebra $A$, the Rees
algebra $R(A)$ is the graded algebra $\oplus_n F^nA$
with the obvious degree preserving multiplication.

The algebra $R(\widehat {\bold H}(q))$ is a formal graded deformation
of $R(H(1,q))$. Let $Z_R'(1,q)$ be the quotient of the center
of this deformation by the ideal generated by $\tau_{kj}$.
Obviously, $Z_R'(1,q)$ is a graded subalgebra
of $Z_R(1,q):=R(Z(1,q))$, and our job is to show that they coincide.

We will use a simple lemma from commutative algebra. 

{\bf Lemma.} Let $A$ be a finitely generated $\Bbb Z_+$-graded
algebra over $\Bbb C$
without zero divisors, and $B\subset A$ be a graded subalgebra.
Suppose that $r$ is a positive integer, and there exists a
constant $C>0$ such that $\dim B[n]>C n^{r-1}$ for large enough
$n$. Then the transcendence degree of $B$ is at least $r$. 

{\bf Proof of the Lemma.} Let $d$ be the transcendence degree of
$B$, and 
$x_1,...,x_d\in B$ be algebraically independent homogeneous
elements.
Let $E=\Bbb C[x_1,...,x_d]$. Let $X=Proj(A)$, and $Y=Proj(E)$
($Y$ is a weighted projective space). 
We have a natural dominant morphism $f: {\rm Cone}(X)\to {\rm
Cone}(Y)$, defined by the embedding $E\subset A$. Any element
$x\in B$ is algebraic over $E$, which means that 
when regarded as a function on ${\rm Cone}(X)$, it is locally
constant 
on a generic fiber of $f$. Let $Z$ be a closed subvariety of  
$X$ of dimension equal to the dimension of $Y$ such that ${\rm
Cone}(Z)$ is transversal 
to a generic fiber of $f$. In this case, the function $x\in B$ on
${\rm Cone}(X)$ is completely determined by its restriction to
${\rm Cone}(Z)$. 
Thus if $R$ is the homogeneous coordinate ring of $Z$ 
then the natural map $B\to R$ is an embedding of graded
algebras. 
This implies that $\dim R[n]>C n^{r-1}$ for large $n$, and hence
$\dim Z\ge r-1$.
But $\dim Z=\dim Y=d-1$, so $d\ge r$, as desired. \qed

Now, the algebras $Z_R(1,q)$ and $Z_R'(1,q)$ are domains. Also,
recall that for generic $t$, $H(t,q)$ is an Azumaya algebra. Thus
$Z_R'(1,q)$  has quadratic growth 
and hence by the lemma has
transcendence degree 3 over $\Bbb C$. This means that
$Z_R(1,q)$ is algebraic over $Z_R'(1,q)$.

Now we will need the following lemma.

\begin{lemma}\label{centdef} Let $A_0$ be an algebra over $\Bbb C$
with center $Z_0$. Assume that $Z_0$ is a domain, and
any derivation of $A_0$ which vanishes on $Z_0$ is inner.
Let $A$ be a flat deformation
of $A_0$ over $\Bbb C[[\hbar]]$. Let $Z$ be the center of $A$ and
$Z_0'=Z/\hbar Z$. Assume that $Z_0$ is algebraic over
$Z_0'$. Then $Z_0=Z_0'$.
\end{lemma}

{\bf Remark.} Note that this lemma is false over a filed $k$ of positive
charateristic. The classical counterexample: $A$ 
is the Weyl algebra with generators $x,y$ and defining
relation $xy-yx=1$, $A_0=k[x,y]$. 

\begin{proof} 
The Hochschild complex of $A[\hbar^{-1}]$ is filtered by degrees in
$\hbar$. There is a Brylinski spectral sequence \cite{Br}
attached to this filtration.
The $E_1$ term of this sequence is the Hochschild
cohomology of $A_0$: $E_1^{p,q}=H^{p+q}(A_0)$.
In particular, $E_1^{p,-p}=Z_0$.
Thus our job is to show that all the differentials
$d_i^{p,-p}$ are zero.

Assume this is not the case, and let $d_n=d_n^{p,-p}$
be the first nonzero differential.
For any $z\in Z_0$, $d_n(z)$ is a coset of derivations of $A_0$ modulo inner
derivations. In particular, $d_n(z)$ gives rise to a well defined derivation
$\tilde d_n(z)$ of $Z_0$.

We claim that $\tilde d_n(z)=0$. Indeed,
$\tilde d_n(z)$ is a derivation with respect to $z$.
If $z\in Z_0'$, then $\tilde d_n(z)=0$. If $z\in Z_0$ and
$P$ is a minimal polynomial
of $z$ over $Z_0'$ then we find $0=\tilde d_n(P(z))(w)=P'(z)
\tilde d_n(z)(w)$.
Since $Z_0$ is a domain, we find $\tilde d_n(z)=0$, as desired.

Thus $d_n(z)$ acts trivially on the center, and hence
by our assumption $d_n(z)=0$. The lemma is proved.
\end{proof}

Let $\tau=\tau(\hbar)$ be a formal path. 
We apply Lemma \ref{centdef} to $A_0=R(H(1,q))$,
and $A=R(\widehat {\bold H}(q))|_{\tau=\tau(\hbar)}$.
The conditions of the lemma hold because
$A_0$ is an Azumaya algebra everywhere
over $Spec Z_R(1,q)$ except a subset of codimension 2,
and any derivation of an Azumaya algebra which is zero on the center is
inner. The conclusion of the lemma implies claim (i).
\end{proof}

Let $u_{31}$ be a simple root of the minimal polynomial of $c$,
and $e$ be the idempotent defined in Theorem \ref{assgrad} (ii).

\begin{corollary}\label{Satake} (``Satake isomorphism'') The natural map
$\phi: Z(t,q)\to eH(t,q)e$ from the center to the spherical subalgebra
given by $z\to ze$ is injective. It defines an isomorphism of 
$Z(t,q)$ onto the center of the spherical subalgebra. 
If $q=1$, $\phi$ is an isomorphism.
\end{corollary}

\begin{proof} The first statement holds because it holds at the
graded level by Theorem \ref{assgrad}.
The second and third statements follow from the first one and the
Poincar\'e series consideration.
\end{proof}

Now fix $q\in \Bbb C^*$ (not necessarily a root of unity), and 
let $t$ be generic in the Zariski sense.  

\begin{proposition}\label{otherfil} 

(i) for $E_6$, ${\rm gr}_{F_{111}}(eH(t,q)e)$
is isomorphic to $B(X,\sigma,{\mathcal L}^{\otimes 3})$.
 
(ii) for $E_7$, ${\rm gr}_{F_{112}}(eH(t,q)e)$
is isomorphic to $B(X,\sigma,{\mathcal L}^{\otimes 2})$. 
\end{proposition}

\begin{proof} The statements are easy in the group case  $t=1$. 
On the other hand, in both cases the family of algebras 
${\rm gr}(eH(t,q)e)$ is flat, since by Theorem \ref{pbw1}, 
their Poincar\'e series is independent of $t,q$.
Therefore, it follows from the theory of noncommutative curves
(see \cite{SV}) that for generic $t$, 
the algebra ${\rm gr}(eH(t,q)e)$ is isomorphic 
to $B(X_{t,q},\sigma_{t,q},{\mathcal L}^{\otimes p})$, where $p=3$ in (i) and $p=2$
in (ii), $X_{t,q}$ is a genus $1$
curve, $\sigma_{t,q}$ its automorphism, and ${\mathcal L}={\mathcal
O}(P)$, where $P$ is a smooth point of $X_{t,q}$. Here
$(X_{t,q},\sigma_{t,q})$ depends algebraically on $t,q$, and
$(X_{1,q},\sigma_{1,q})=(X,q^\ell)$, as explained above. 

If $q$ is a root of unity, such that $q^\ell$ has degree $N$, 
then $H(t,q)$ is an Azumaya algebra of rank $\ell N$,
and hence $\sigma_{t,q}$ has order $N$. 

Let $\Sigma_N$ be the modular curve
parametrizing elliptic curves with points of order 
$N$. We see that if $q$ is a root of unity as above,
we get a regular map $t\mapsto \psi_q(t)$ from an open subset of
$\Bbb T$ to $\Sigma_N$, given by $\psi_q(t)=(X_{t,q},\sigma_{t,q})$.
For large enough $N$, the curve $\Sigma_N$ is not rational, and
hence the map $\psi_q$ must be constant. Thus, $X_{t,q}=X$ for
all $t,q$. So we can think of $\sigma_{t,q}$ as an element of
$\Bbb C^*$. 

If $q$ is a root of unity, so must be $\sigma_{t,q}$. Therefore,
for $q$ being a root of unity, $\sigma_{t,q}$ is independent of $t$. 
Thus, for such $q$, $\sigma_{t,q}=\sigma_{1,q}=q^\ell$. 
Since roots of unity are Zariski dense, this equality holds for
all $q$. We are done. 
\end{proof}

This allows one to give the following descriptions 
of the associated graded algebras of $Z(t,q)$ for Zariski generic
$t$.

\begin{corollary}\label{otherfilcen}
(i) for $D_4$, ${\rm gr}_{F_{111}}(Z(t,q))$
is isomorphic to $B(X_3,1,{\mathcal L})$,
with grading multiplied by $N$.
So the Poincar\'e series of $Z(t,q)$ under this filtration is 
$$
1+\frac{3z^N}{(1-z^N)^2}.
$$

(ii) for $E_6$, ${\rm gr}_{F_{111}}(Z(t,q))$
is isomorphic to $B(X,1,{\mathcal L}^{\otimes 3})$,
with grading multiplied by $N$.
So the Poincar\'e series of $Z(t,q)$ under this filtration is 
$$
1+\frac{3z^N}{(1-z^N)^2}.
$$

(iii) for $E_7$, ${\rm gr}_{F_{112}}(Z(t,q))$
is isomorphic to $B(X,1,{\mathcal L}^{\otimes 2})$,
with grading multiplied by $N$.
So the Poincar\'e series of $Z(t,q)$ under this filtration is 
$$
1+\frac{2z^N}{(1-z^N)^2}.
$$

(iv) for $E_8$, ${\rm gr}_{F_{123}}(Z(t,q))$
is isomorphic to $B(X,1,{\mathcal L})$,
with grading multiplied by $N$.
So the Poincar\'e series of $Z(t,q)$ under this filtration is 
$$
1+\frac{z^N}{(1-z^N)^2}.
$$
\end{corollary}

\begin{proof} Statements (i) and (iv) follow from Theorem \ref{assgrad} 
and Theorem \ref{grcenter}. Statements (ii) and (iii) 
follow from Proposition \ref{otherfil} and
Corollary \ref{Satake}.
\end{proof}

\subsection{Quantization of del Pezzo surfaces}

Recall that a del Pezzo surface is a smooth projective surface
with ample anticanonical bundle $K^{-1}$. Apart from $\Bbb P^1\times \Bbb P^1$
such surfaces are obtained from $\Bbb P^2$ by blowing up
$n$ sufficiently generic points ($n\le 8$).
The degree of a del Pezzo surface is the dimension of
the projective space of $\Gamma(K^{-1})$.
For example, a del Pezzo surface of degree 3 is a cubic surface
in $\Bbb P^3$. The degree of a projective plane
with $n$ generic points blown up is $9-n$.

Let $q$ be a root of unity, such that the order of $q^\ell$ is
$N$. Let $t$ be generic in the Zariski sense. 

\begin{theorem}\label{thm:Ztqgen}
(i) In the $D_4$ case, $Z(t,q)$
is generated by degree $N$ elements $x,y,z$
with defining relation
$$
xyz+x^2+y^2+z^2+a_2x+a_3y+a_4z+a_1=0,
$$
where $a_m$ are functions of $t_{kj}$;

(ii) In the $E_6$ case, $Z(t,q)$
is generated by degree $N$ elements $x,y,z$
with defining relation
$$
xyz+x^3+y^3+z^2+a_1x^2+a_2y^2+a_3x+a_4y+a_5z+a_6=0,
$$
where $a_m$ are functions of $t_{kj}$;

(iii) In the $E_7$ case, $Z(t,q)$
is generated by degree $N$ elements $x,y$ and degree $2N$
element $z$
with defining relation
$$
xyz+x^4+y^3+z^2+a_1x^3+a_2x^2+a_3y^2+a_4x+a_5y+a_6z+a_7=0,
$$
where $a_m$ are functions of $t_{kj}$;

(iv) In the $E_8$ case, $Z(t,q)$
is generated by degree $N$ element $x$, degree
2N element $y$, and degree $3N$
element $z$ with defining relation
$$
xyz+x^5+y^3+z^2+a_1x^4+a_2y^2+a_3x^3+a_4z+a_5x^2+a_6y+a_7x+a_8=0,
$$
where $a_m$ are functions of $t_{kj}$.

Here the degree in cases (i), (ii), (iii), (iv) is computed with respect 
to the filtration $F_{111}$, $F_{111}$, $F_{112}$, $F_{123}$,
respectively. 
\end{theorem}

\begin{proof} The proof is based on Corollary \ref{otherfilcen}.

Namely, in the $D_4$ case,
$Z(t,q)$ is generated by
three elements $x,y,z$ of degree $N$
(since this is true already for the graded algebra, by Theorem 
\ref{assgrad}).
 From looking at the Poincar\'e series
it is clear that these elements satisfy a cubic defining relation,
whose homogeneous part (as we know from studying the graded algebras)
is $xyz=0$ (triangle at infinity).
Thus, this relation
has the form
$$
xyz+Q(x,y,z)=0,
$$
where $Q$ is an inhomogeneous quadratic form.
Using affine transformations (i.e., shifts of $x,y,z$ by constants),
this equation can be brought to the form
$$
xyz+A(x)+B(y)+C(z)=0
$$
where $A,B,C$ are (at most) quadratic polynomials.
For generic $t$ these polynomials
have nonvanishing quadratic terms, since 
this is so for $t=1$. Thus the coefficients 
in these terms can be normalized to be $1$, which gives the
result. 

In the $E_6$ case,
$Z(t,q)$ is also generated by
three elements $x,y,z$ of degree $N$
(since this is true already for the graded algebra).
>From looking at the Poincar\'e series
it is clear that these elements satisfy a cubic defining relation.
However, now the curve at infinity is a nodal $\Bbb P^1$ rather than
triangle, so the homogeneous part of the cubic relation
can be brought by linear transformations
to the form $xyz+x^3+y^3=0$.
Thus, the cubic relation
has the form
$$
xyz+x^3+y^3+Q(x,y,z)=0,
$$
where $Q$ is an inhomogeneous quadratic form.
Using affine transformations,
this equation can be brought to the form
$$
xyz+A(x)+B(y)+C(z)=0
$$
where $A,B$ are cubic polynomials and $C$ is at most quadratic.
For generic $t$, the polynomial $C$
has nonvanishing quadratic coefficient, since it is so for
$t=1$. So the leading coefficients of $A,B,C$ can be normalized
to be $1$. 

In the $E_7$ case,
$Z(t,q)$ is generated by
three elements $x,y,z$ of degree $N, N, 2N$
(since this is true already for the graded algebra).
>From looking at the Poincar\'e series
it is clear that these elements satisfy a defining relation
in degree 4. Since the curve at infinity is a nodal $\Bbb P^1$, the
homogeneous part of this relation
can be brought by linear transformations
to the form $xyz+x^4+z^2=0$.
Thus, after linear transformation the inhomogeneous relation
can be brought to the form
$$
xyz+A(x)+B(y)+C(z)=0
$$
where $A$ is a quartic polynomial,
$C$ a quadratic polynomial, and
$B$ is at most cubic.
For generic $t$ the polynomial $B$
has nonvanishing cubic term, since it is so for $t=1$. 
So the leading coefficients of $A,B,C$ can be normalized
to be $1$. 

In the $E_8$ case,
$Z(t,q)$ is generated by
three elements $x,y,z$ of degree $N, 2N, 3N$
(since this is true already for the graded algebra).
>From looking at the Poincar\'e series
it is clear that these elements satisfy a defining relation
in degree 6. Since the curve at infinity is a nodal $\Bbb P^1$, the
homogeneous part of this relation
can be brought by linear transformations
to the form $xyz+y^3+z^2=0$.
Thus, after linear transformation the inhomogeneous relation
can be brought to the form
$$
xyz+A(x)+B(y)+C(z)=0
$$
where $A$ is at most a quintic polynomial,
$C$ a quadratic polynomial, and
$B$ is a cubic polynomial.
For generic $t$ the polynomial $A$
has nonvanishing quintic term, since it is so
for $t=1$. So the leading coefficients of $A,B,C$ can be normalized
to be $1$. 
\end{proof}

The theorem shows that we can view $S(t,q)$
of types $D_4$, $E_6$, $E_7$, $E_8$ as sitting inside
the weighted projective space with weights
(1,1,1),(1,1,1),(1,1,2),(1,2,3), respectively.
Denote by $\overline{S(t,q)}$ its closure
(=Proj of the Rees algebra of $Z(t,q)$ for the filtrations 
$F_{111},F_{111},F_{112},F_{123}$, respectively),
and by $C(t,q)$ the curve at infinity
(=Proj of the graded algebra of $Z(t,q)$
with respect to these filtrations).

\begin{corollary}
In the cases (i)-(iv), the curve at infinity $C(t,q)$ consists
of smooth points of the compact surface $\overline{S(t,q)}$.
Therefore, the surface $\overline{S(t,q)}$ is smooth for
generic $t$.  
\end{corollary}

\begin{proof}
This follows from the nonvanishing of the highest coefficients
of $A,B,C$.
\end{proof}

Thus we get the following result on the structure of $S(t,q)$ 
for generic $t$.

\begin{theorem}\label{dp} $S(t,q)$ is isomorphic to:

(i) \cite{Ob} in the $D_4$ case, a Del Pezzo surface
of degree 3 with a triangle  removed;

(ii) in the $E_6$ case, a Del Pezzo surface
of degree 3 with a nodal $\Bbb P^1$  removed;

(iii) in the $E_7$ case, a Del Pezzo surface
of degree 2 with a nodal $\Bbb P^1$  removed;

(iv) in the $E_8$ case, a Del Pezzo surface
of degree 1 with a nodal $\Bbb P^1$  removed.
\end{theorem}

\begin{proof}
The theorem follows from the well known fact in algebraic geometry
(see e.g. \cite{D}, p. 60-71) that equations of the form $xyz+A(x)+B(y)+C(z)$
in weighted projective space give realizations of Del Pezzo
surfaces. Namely, such realizations of Del Pezzo surfaces $S$ are 
obtained by considering the Proj of the ring $\oplus_{k\ge 0}H^0(S,(K^*)^{\otimes
k})$, where $K^*$ is the anticanonical bundle on $S$. 
\end{proof}

Abusing terminology, we will use the name ``del Pezzo surface''
for all, and not only smooth, members of the family $S(t,q)$.

Note that any smooth surface $S$ in $\Bbb C^3$ has a natural
symplectic structure up to scaling. If the equation of the surface is
$F(x,y,z)=0$ then the symplectic form is $\omega=\frac{dx\wedge
dy\wedge dz}{dF}$. If $S$ is singular, this symplectic form 
becomes singular at singular points of $S$, but still defines a
regular Poisson structure on $S$. 
It is defined by the formulas 
$$
\lbrace{ x,y\rbrace}=\frac{\partial F}{\partial z},
\lbrace{ y,z\rbrace}=\frac{\partial F}{\partial x},
\lbrace{ z,x\rbrace}=\frac{\partial F}{\partial y}.
$$

Theorem \ref{dp} and Corollary \ref{Satake} implies

\begin{theorem} \label{quant} 
Let $q=e^{\hbar}$, where $\hbar$ is a formal parameter.
Then the algebra $eH(t,q)e$ is a deformation quantization of
the del Pezzo surface $S(t,1)$, with its natural Poisson
structure (with an appropriate normalization).
\end{theorem}

\begin{proof}
It is sufficient to prove the result for generic $t$. 
In this case the surface $S(t,1)$ is smooth and the natural Poisson
structure is symplectic, given by the formulas above.

Now, the algebra $eH(t,q)e$ is a quantization of some 
(maybe different) Poisson
structure $\lbrace{,\rbrace}'$ on $S(t,1)$. This Poisson
structure must have the form
$\lbrace{,\rbrace}'=f\lbrace{,\rbrace}$,
where $f$ is a polynomial function on $S(t,1)$. 
But $\lbrace{,\rbrace}'$ must preserve the filtration
on $eH(t,1)e$, which implies that $f$ has to be constant. 
The theorem is proved. 
\end{proof}

\begin{proposition}\label{simcon}
For generic $t$, $H^1(S(t,1),\Bbb C)=0$, 
and $H^2(S(t,1),\Bbb C)=\Bbb C^{r+1}$, where $r$ is the rank of the
Dynkin diagram $D$ (i.e., $4,6,7,8$). 
\end{proposition}

\begin{proof} The identity $H^1(S(t,1),\Bbb C)=0$ obviously holds
for $t=1$, since in this case $S(t,1)$ is $T/\Bbb Z_\ell$. 
Thus the only way $H^1(S(t,1),\Bbb C)$ could be nonzero 
for generic $t\ne 1$ is if there were vanishing 1-cycles as $t\to 1$
near the singular points of $S(1,1)$ (this follows from the fact
that the deformation $S(t,1)$ near $t=1$ is topologically trivial
everywhere including infinity, except the singular points of
$T/\Bbb Z_\ell$). But it is clear that at the singular points,
there could only be vanishing 2-cycles and not 1-cycles. 
Thus, $H^1(S(t,1),\Bbb C)=0$ for generic $t$.

Now, it is easy to see that the Euler characteristic of $S(t,1)$
is $r+2$ for generic $t$. This implies the result, since $S(t,1)$ is affine and
can't have cohomology above degree 2.  
\end{proof}

{\bf Remark 1.} This proposition shows that the symplectic
structure on $S(t,1)$ for generic $t$ is
unique up to scaling, since by Proposition \ref{simcon}, 
$H^1(S(t,1),\Bbb C)=0$ and hence $S(t,1)$ does not have 
nonconstant nowhere vanishing functions
(the logarithmic differential of such a function would represent a
nontrivial class in $H^1$). 

{\bf Remark 2.} Since our algebras are equipped with filtrations,
we also get quantum surfaces in the sense of noncommutative algebraic
geometry, which are quantum deformations of commutative compact
surfaces. Namely, the homogeneous coordinate rings of these quantum surfaces
are the Rees algebras of $eH(t,q)e$ equipped with filtrations $F_{ijk}$.

Let us specify which commutative surfaces are quantized in this
way. Obviously, these commutative surfaces are the Projs of the Rees algebras
of $eH(t,1)e=Z(t,1)$. So let us describe (omitting the proofs)
what surfaces we get (for generic $t$).
The cases of $F_{111}$ for $D_4$, $F_{111}$ for $E_6$, $F_{112}$
for $E_7$, and $F_{123}$ for $E_8$ are covered by Theorem
\ref{dp}. The remaining cases are $F_{112}$ for $D_4$, $F_{112}$ and
$F_{123}$ for $E_6$ and $F_{123}$ for $E_7$.

$\bullet$ $F_{112}$ for $D_4$: 
$\overline{S(t,1)}$ is a singular del Pezzo surface of degree 2, 
the divisor at infinity
consists of two rational curves intersecting at two points,
and both intersection points carry $A_1$ singularities of the
surface. 

$\bullet$ $F_{112}$ for $E_6$: $\overline{S(t,1)}$ is
a smooth Del Pezzo surface of
degree $2$ (i.e., of type $E_7$), the divisor at infinity
consists of two rational curves intersecting at two points. 

$\bullet$ $F_{123}$ for $E_6$: $\overline{S(t,1)}$ is
a singular del Pezzo surface of degree $1$ (i.e., of type $E_8$), 
the divisor at infinity is a rational 
curve intersecting itself at a point.
The surface has a singularity of type $A_2$ at this point.  
   
$\bullet$ $F_{123}$ for $E_7$: 
$\overline{S(t,1)}$ is
a singular del Pezzo surface of degree $1$ (i.e., of type $E_8$), 
the divisor at infinity is a rational 
curve intersecting itself at a point.
The surface has a singularity of type $A_1$ at this point.  

Finally, we derive a corollary about the Hochschild cohomology
of the generalized Cherednik algebras. 

\begin{proposition}
Let $q=e^\hbar$, where $\hbar$ is a formal parameter. 
Then the Betti numbers of the Hochschild cohomology of the algebra
$H(t,q)[\hbar^{-1}]$ for generic $t$ are $b_0=1$, $b_1=0$,
$b_2=r+1$, $b_i=0$ for $i>2$. 
\end{proposition}

\begin{proof} The proof is analogous to the proof in the case of
$D_4$, given in \cite{Ob}. Namely, 
the result follows from Proposition \ref{simcon}, Theorem
\ref{quant}, and the theorem 
on cohomology of quantizations of symplectic manifolds (see e.g. 
\cite{Do} and references therein). 
\end{proof}

\subsection{Nondegeneracy of the map $t\to S(t,q)$}

Del Pezzo surfaces given by the equations $xyz+A(x)+B(y)+C(z)$
form a moduli space ${\mathcal M}$, coordinates on which
are the coefficients of $A,B,C$ (stipulating that
the highest coefficients are $1$
and there is only one independent coefficient among the constant terms of
$A,B,C$). Thus the dimension of
${\mathcal M}$ is equal to the rank $r$ of the corresponding Dynkin diagram
$D$ (i.e., 4,6,7,8).

Let $q$ be a root of unity, and
$\zeta_q: {\mathcal U}\to {\mathcal M}$ be the map
attaching $S(t,q)$ to $t$ (here ${\mathcal U}$ is some open set
in $\Bbb T$).

\begin{theorem}\label{domi}
The map $\zeta_q$ is dominant.
\end{theorem}

\begin{proof}
Consider a formal path $t=t(s) $ such that $t(0)=1$
and $S(t(s),q)=S(1,q)$. Thus
$H(t(s),q)$ is a formal deformation of $H(1,q)$ over
a fixed center $Z(1,q)$. We will show that this deformation is trivial.

Consider the first nontrivial order of the $s$-expansion.
In this order our deformation defines a Hochschild 2-cocycle
$\gamma$ of $H(1,q)$ as an algebra over $Z(1,q)$.
It suffices to show that this cocycle is trivial; then we can
make a gauge transformation to insure that the lowest nontrivial order in
$s$ is becomes one step higher, and obtain the result by induction.

The group $H^2_{Z(1,q)}(H(1,q))$ is a finitely generated
module over $Z(1,q)$,
i.e. a coherent sheaf on $S(1,q)$.
The surface $S(1,q)$ is isomorphic to $T/\Bbb Z_\ell$, and hence
has $m$ isolated singular points ($m=4,3,3,3$).
Near any other point, $H(1,q)$ (and $H(t,q)$ for $t$ close to $1$)
is an Azumaya algebra.
This shows that $\gamma$ vanishes outside
of the singular points. Hence in showing that $\gamma=0$, we
may replace $H(1,q)$ and $Z(1,q)$ by their completions
near the singular points.

Now, let $p$ be a singular point of $S(1,q)$, and $\Bbb Z_k$
be the stabilizer of this point in $\Bbb Z_\ell$.
It is clear that the completion of $H(1,q)$ near $p$
is Morita equivalent to $A=\Bbb C[\Bbb Z_k]\ltimes \Bbb C[[x,y]]$,
where the generator of $\Bbb Z_k$ multiplies $x$ by
the primitive $k$-th root of unity, and $y$ by the inverse of this root.
The cocycle $\gamma|_A$ comes from a formal deformation $A(s)$ of $A$ that
keeps the center fixed. 

Now, it is easy to compute using Koszul resolutions that
$$
H^2(\Bbb C[\Bbb Z_k]\ltimes \Bbb C[[x,y]])=
H^2(\Bbb C[[x,y]])^{\Bbb Z_k}\oplus \Bbb C^{k-1}.
$$
If $\gamma$ had a nontrivial projection to the first summand,
then the center would collapse under deformation. 
Thus, $\gamma$ belongs to the second summand. 
This means that our deformation falls into the family 
of algebras from \cite{CBH}. For this family, it is known that 
the deformation of the center is a versal deformation of
the singularity $A_{k-1}$. This means that if the center is fixed
then the deformation is trivial. Thus we see 
that $A(s)$ is a trivial deformation, and hence $\gamma|_A=0$.
Thus $\gamma=0$.

Now Lemma \ref{nontriv} implies that the path $t=t(s)$ must be trivial:
$t(s)=1$. This implies the statement.
\end{proof}

\subsection{A linear algebra application} 

The results of Section 6 imply the following 
result from linear algebra, which appears to be new. 
Fix a diagram $D$ of type $D_4$, $E_6$, $E_7$, $E_8$. 

\begin{theorem} Let $q$ be a root of unity such that $q^\ell$ has
order $N$, and $t$ be generic (more specifically, $t\in U$ 
with $u_{kj}\ne u_{kj'}$ for $j\ne j'$). Let $S(t,q)$ be the space of 
conjugacy classes of collections of diagonalizable matrices 
$T_1,...,T_m$ of size $\ell N$, such that eigenvalues of $T_k$
are $u_{kj}$, $j=1,...,d_k$ (with some multiplicities), and
$T_1...T_m=q$. Then the multiplicities of the eigenvalues are all
the same (i.e., equal $\ell N/d_k$), and $S(t,q)$ is an affine Del Pezzo surface
described in Theorem \ref{dp}. Moreover, a generic surface 
of this kind is obtained in this way.   
\end{theorem} 

\section{The Riemann-Hilbert map}

\subsection{Preprojective algebras}

Let $Q$ be a quiver, and let $E(Q)$ be the set of edges of $Q$. 
The path algebra of $Q$ is spanned by paths in $Q$
with multiplication given by concatenation of paths.
In particular, it contains the idempotents $p_i$ corresponding
to the paths of length $0$ at the vertices $i$ of $Q$.

Let ${\Bbb D}(Q)$ be the double of $Q$, obtained by
adding, for any edge $h\in E(Q)$, a new edge $h^*$ in the opposite
direction. Setting $h^{**}=h$, we get an involution of the set
of edges $E({\Bbb D}(Q))$. 

The Gelfand-Ponomarev deformed preprojective algebra
$\Pi_\mu$ is the quotient of the path algebra of ${\Bbb D}(Q)$ by the relation
$$
\sum_{h\in E(Q)} [h,h^*]=\sum_i \mu_ip_i, 
$$
where $\mu_i$ are complex numbers.

Let $\Gamma$ be a finite subgroup of $SL(2,\Bbb C)$.
To such a group, Crawley-Boevey and Holland \cite{CBH} assigned
an algebra $Q_\Gamma(c)$ generated by
$\Gamma$ and its tautological 2-dimensional representation
$V$ with defining relations
$$
gvg^{-1}=v^g, g\in \Gamma, v\in V,
$$
and
$$
[v,w]=(v,w)\sum_g c_gg,
$$
where $c$ is a class function on $\Gamma$ and
$(v,w)$ is the symplectic inner product. Let $e_i$ be primitive idempotents
in $\Bbb C\Gamma$ attached to irreducible representations
$V_i$ of $\Gamma$ (they are unique up to conjugation).
We denote by $Q_\Gamma'(c)$ the algebra $\oplus_{i,j}e_iQ_\Gamma(c) e_j$.

Let $Q$ be an affine quiver.
McKay's correspondence assigns to $Q$ a finite subgroup $\Gamma$ of $SL(2)$,
whose irreducible representations are labeled by vertices of $Q$.
Let $\chi_i$ be the character of the $i$-th irreducible representation.

\begin{proposition} \cite{CBH}
The algebra $\Pi_\mu$ is isomorphic
to $Q_\Gamma'(c)$ if $\sum c_g\chi_i(g)=\mu_i$.
\end{proposition}

Now assume that $Q$ is star-like (i.e, $D_4,E_6,E_7,E_8$).
In this case, let $i_0$ be the nodal vertex, and
$p=p_{i_0}$ be the correspondent idempotent of $\Pi_\mu$.
Let $K(\mu)=p\Pi_\mu p$.

\begin{proposition}\label{mov} (see \cite{MOV,Me,CB1})
The algebra $K(\mu)$ is generated by elements $U_k$, $k=1,...,m$,
corresponding to the legs of $Q$, subject to defining relations
$$
U_k(U_k-\mu_{i_1(k)})(U_k-\mu_{i_1(k)}-...-\mu_{i_{d_k-1}(k)})=0,
$$
where $i_1(k),...,i_{d_k-1}(k)$ are the vertices
of the $k$-th leg of $Q$ enumerated from the nodal vertex, and
$$
\sum_{k=1}^m U_k=-\mu_{i_0}.
$$
\end{proposition}

\begin{proof} The elements $U_k$ are just the elements
$h_k^*h_k$, where $h_k$ are the edges of ${\Bbb D}(Q)$ starting
at $i_0$ and going
along the $k$-th leg. They obviously generate $K(\mu)$.
It is easy to compute that $U_k$ satisfy
the relations above, and it is not hard to check that
these relations are defining.
\end{proof}

Thus, the algebra $K(\mu)$ is an additive analog of
the algebra $H(t,q)$.

{\bf Remark.} A recent paper \cite{CBS} introduces a multiplicative
analog of the preprojective algebra -- the multiplicative
preprojective algebra $\Pi^{\rm mult}_\mu$ of a quiver $Q$. Like the usual
preprojective algebra, this algebra has 
idempotents $p_i$ attached to the vertices $i$ of the quiver. 
It can be shown (see the appendix below) that the algebra $H(t,q)$ is
isomorphic to $p\Pi_\lambda^{\rm mult}p$ (for appropriate
$\lambda$). 

\subsection{The Riemann-Hilbert map}

Let $\bold K$ be the algebra $K(\mu)$ where $\mu$ are formal parameters
(i.e. it is an algebra over $\Bbb C[[\mu]]$).
Let $\bold H$ be the algebra $H(t,q)$ with
$q=e^\hbar$ and $u_{kj}=e^{\beta_{kj}}$ where
$\beta_{kj}$ are formal parameters.
Note that $\bold H$ is different from the algebras considered
in Sections 2-6, since now we take completion near the point
$u=1$ (unipotent case) rather than $t=1$ (infinite group case). 

Representations of the algebra $\bold K$ are solutions
of the additive Deligne-Simpson problem, while
representations of $\bold H$ are solutions of the multiplicative one.
Thus we have a Riemann-Hilbert map between 
completions of these algebras, defined as follows.

Let $z_1,...,z_m$ be distinct points on $\Bbb C \Bbb P_1$.
Consider the flat connection $\nabla$
on the trivial bundle over $X=\Bbb P_1\setminus
\lbrace{z_1,...,z_m\rbrace}$ with fiber $\bold K$
which has
first order poles with residues $U_k+\mu_{i_0}/m$ at $z_k$.
Let $z_0\ne z_k$ for any $k$ and $\gamma_k$
be the standard generators of $\pi_1(X,z_0)$ going around $z_k$
such that $\prod_k \gamma_k=1$. Let
$\widetilde{\bold K}$ be the completion of $\bold K$ with respect to the ideal
defined by $U_1,...,U_m$, and
$T_k'\in \widetilde{\bold K}$ be the monodromies of the
connection $\nabla$ around $\gamma_k$.
Let $\bar T_k=e^{2\pi i\lambda_k}T_k'$, where
$$
\lambda_k=-\mu_{i_0}/m-
\sum_{j=1}^{d_k-1}\frac{j}{d_k}\mu_{i_{d_k-j}(k)}.
$$
Denote by $e^{\bar \beta_{kj}}$ the roots of the polynomial
equation satisfied by $\bar T_k$, and by
$e^{\bar \hbar}$ the product of $T_k$. These are
exponentials of linear functions of $\mu$.

Let $\widetilde{\bold H}$ be the completion of ${\bold H}$ with respect to
the ideal generated by the elements $T_k-1$.

\begin{proposition}
The map $T_k\to \bar T_k$,
$\beta_{kj}\to \bar \beta_{kj}$, $\hbar\to \bar\hbar$
is an embedding $\phi_0: {\bold H}\to \widetilde{\bold K}$
which extends by continuity to an isomorphism
$\phi: \widetilde{\bold H}\to \widetilde{\bold K}$.
\end{proposition}

\begin{proof}
The fact that the given formulas define 
a homomorphism of algebras is obtained by an easy direct
computation. The fact that $\phi_0$ is injective 
may be checked for $\hbar=0$ (i.e., modulo $\hbar$), since the target algebra 
$\widetilde{\bold K}$ is flat over $\Bbb C[[\hbar]]$. If
$\hbar=0$, then the injectivity of $\phi_0$ follows from the fact that 
the spectrum $S(t,1)$ of the center of $H(t,1)$ is irreducible 
(Proposition \ref{irraf}), and therefore the center embeds into
its completion at every point. The fact that $\phi_0$ extends to an
isomorphism of completions is now straightforward.  
\end{proof} 

{\bf Remark.}
The isomorphism $\phi$ is independent on the choice of the point $z_0$ up to
inner automorphisms. Thus the corresponding isomorphism 
between the centers $\phi_c: Z(\widetilde{\bold H})\to Z(\widetilde{\bold
K})$ is independent on the choice of $z_0$. 
As a function of $z_1,...,z_m$, $\phi_c$ 
is projectively invariant. Therefore, 
it is completely canonical in the cases $E_6,E_7,E_8$ and 
depends on one parameter (cross ratio) for $D_4$.

For numerical values of $\mu$ and $t,q$, the map $\phi$ is not well defined,
since we cannot compute monodromies of connections with infinite dimensional
fiber. However, we have the following proposition.

\begin{proposition}
For parameters related as above, we have
a ``pullback'' functor between the categories of finite
dimensional representations \linebreak $\phi^*: {\rm Rep}_f(K(\mu))\to
{\rm Rep}_f H(t,q)$.
\end{proposition}

This functor is obviously far from being an equivalence, since the
2-parameter family of finite dimensional representations
at $q$ being a root of unity (with $q^\ell\ne 1$) does not belong to its image.
However, the restriction of $\phi^*$ to the situation
when $c_1=0$ and $q=1$ produces a canonical
holomorphic (but not algebraic) mapping $\phi^*: S_0(\mu)\to S(t,1)$,
where $S_0(\mu)$ is the versal deformation of $\Bbb C^2/\Gamma$
(=the spectrum of the center of $K(\mu)$); this map is an isomorphism
near the origin. More precisely, as we mentioned above, in the case of
$D_4$ the map $\phi^*$ depends on the cross
ratio $s$ of the points $z_1,..,z_4$, while in types
$E_6,E_7,E_8$ it is independent of any choices
and completely canonical.

{\bf Remark 1.} In the case of $D_4$, consider the inverse map
$(\phi^*)^{-1}$, and the point
$g(s)=(\phi^*)^{-1}(y,t,s)\in S_0(\mu)$, where $y\in S(t,1)$
(here $\mu$ is a linear transformation of $\log t$ as explained above).
Clearly, $g(s)$, regarded as a function of the cross-ratio $s$, represents
an isomonodromic deformation of 2-dimensional local systems on
$\Bbb CP^1$
wit 4 regular singular points. Thus, in appropriate coordinates it
satisfies the differential equation Painlev\'e VI (with general parameters),
and generic solutions of Painlev\'e VI are obtained in this way
(see \cite{Ok}).
The known symmetry of Painlev\'e VI under the affine Weyl group
$\tilde W(D_4)$ of type $D_4$ (\cite{Ok})
is combined from the usual $W(D_4)$ symmetry on the deformation of
$\Bbb C^2/\Gamma$ and the lattice of rank 4
coming from the map $t\to \log t$. 

{\bf Remark 2.} As above, let $\Gamma$ be the group attached to the diagram
$\widehat D$ via McKay's correspondence. Let $\Gamma_+$ be the
quotient of $\Gamma$ by $\pm 1$. 
Then $\Gamma_+$ has generators $a,b,c$ 
and the following defining relations:

$D_4$:  $a^2=b^2=c^2=1$, $abc=1$. 

$E_6$: $a^3=b^3=c^2=1$, $abc=1$. 

$E_7$: $a^2=b^4=c^3=1$, $abc=1$.

$E_8$: $a^2=b^3=c^5=1$, $abc=1$. 

For $D_4$ it is convenient to add a new generator $d$ 
and write the relations as 
$a^2=b^2=c^2=d^1=1$, $abcd=1$.

Then in all cases 
the relations are exactly the same as for the group $G$, except that
the order of 
the last generator equals its order in
$G$ minus 1. 

This shows that any irreducible representation $V$  
of $\Gamma$ can be viewed as a representation
of the algebra $H(t,q)$ for appropriate 
$t$ and $q$. Indeed, the above implies that 
the action of $\Gamma$ in $V$ is generated 
by operators $a,b,c$ satisfying the above relations 
up to sign. But then (rescaled versions of) $a,b,c$ define 
an action in $V$ of a 1-parameter family of algebras $H(t,q)$.
Indeed, the relations for $H(t,q)$ are obtained 
if for $E_{6,7,8}$ one replaces the relation
$c^p=1$ ($p=2,3,5$) by the relation $(c^p-1)(c-\lambda)=0$,
and for $D_4$ one replaces the relation $d=1$ with the
relation $(d-1)(d-\lambda)=0$.

\section{Proofs of Theorems \ref{pbw1} and \ref{assgrad}}

\subsection{Efficient presentations}\label{pbws1}
The main difficulty in doing computations on noncommutative algebras given
a presentation is that there are no general algorithms for such
computations; indeed, most natural problems about finitely presented
algebras are known to be undecidable.  In the commutative case, the primary
tool is the notion of Gr\"obner basis; while this notion has been extended
to the noncommutative case \cite{Gr}, noncommutative Gr\"obner bases
need not be finite in general, and their finiteness depends highly on the
choice of presentation.  Thus the first task in computing in our algebras
is to find ``efficient'' presentations, i.e., presentations admitting
finite Gr\"obner bases.

For computational purposes, it turns out to be convenient to first rescale
the generators slightly.  We thus obtain the following presentations of our
algebras.  In each case, we observe that $a$ can be expressed easily as a
polynomial in $a^{-1} = bc/q$ ($bcd/q$ for $D_4$), so it suffices to take
$b$ and $c$ (and $d$) as generators.  This observation, plus some mild
rescaling of the generators and parameters, gives rise to the following
(not yet efficient) versions of the algebras.

For $D_4$:
\[
\langle b,c,d\mid
c^2-g_1 c+1,
b^2-f_1 b+1,
d^2-h_1 d+1,
(bcd)^2-e_1(bcd)+Q
\rangle,
\]
where $Q=q^2$. (Note that $b$ and $c$ have been rescaled to change the
last coefficients of their minimal polynomials.)

For $E_6$:
\[
\langle b,c\mid
c^3-g_1 c^2+g_2 c-1,
b^3-f_1 b^2+f_2 b-1,
(bc)^3-e_2 (bc)^2+e_1 (bc)-Q
\rangle,
\]
where $Q=q^3$.

For $E_7$:
\[
\langle b,c\mid
c^4-g_1 c^3+g_2 c^2-g_3 c+1,
b^4-f_1 b^3+f_2 b^2-f_3 b+Q,
(bc)^2-e_1 (bc)+Q
\rangle,
\]
where $Q=q^4$.  

For $E_8$:
\[
\langle b,c\mid
c^6-g_1 c^5+g_2 c^4-g_3 c^3+g_4 c^2-g_5 c + 1,
b^3-f_1 b^2 + f_2 b + Q,
(bc)^2 - e_1 (bc) + Q
\rangle,
\]
where $Q=q^6$.  (Again, $b$ and $c$ have been rescaled for $E_6$
and $E_8$.)

For $E_6$, the above filtration already admits a finite Gr\"obner basis,
with respect to the ``shortlex'' term order.  We first note that
$$
\begin{aligned}
b^3 &= f_1 b^2 - f_2 b + 1\\
c^3 &= g_1 c^2 - g_2 c + 1.
\end{aligned}
$$
This allows us to compute inverses in the algebra; multiplying the third
relation by $(bc)^{-1}$ gives
\[
bcbc = e_2 bc - e_1 + Q (c^2 - g_1 c + g_2) (b^2 - f_1 b + f_2);
\]
conjugating by $b$ gives
\[
cbcb = e_2 cb - e_1 + Q (b^2 - f_1 b + f_2) (c^2 - g_1 c + g_2),
\]
allowing us to expand $b^2c^2$ in smaller monomials.  The noncommutative
analogue of the Buchberger algorithm shows that these four relations
generate a Gr\"obner basis, with leading terms $b^3$, $c^3$, $bcbc$,
$b^2c^2$.  In particular, any element in the algebra can be expressed as a
linear combination of words for which none of these four words is a
subword.  Since the elements Gr\"obner basis are compatible with the
filtration (the degree of the leading term is at least as high as the
degree of the remaining terms), and the coefficients are in $\CC[\Bbb T]$, we
conclude that Theorem \ref{pbw1} holds for the 123 filtration of the $E_6$
algebra.  A similar short computation gives a Gr\"obner basis for the above
presentation of $E_7$ compatible with the 123 filtration, proving Theorem
\ref{pbw1} for that case as well.  (In that case, the leading terms of the
Gr\"obner basis are $c^4$, $b^3$, $bcb$, $bccbc$, $bbccc$.)

For $E_8$, it appears that the above presentation does not admit a finite
Gr\"obner basis for any choice of term order compatible with the
filtration.  We must therefore use a different generating set.  For $0\le
i\le 5$, we let $d_i = b c^{-i}$, and take as generators $c$ together with
$d_0$ through $d_5$.  We then obtain (via a much longer computation,
although still quite short on a computer) a shortlex Gr\"obner basis, with
respect to the variable ordering $c,d_1,d_3,d_5,d_0,d_2,d_4$, having
leading terms as follows:
\begin{itemize}
\item[(1)] $c^6$.
\item[(2)] $d_i c$, $0\le i\le 5$
\item[(3)] $d_i d_j$ with $2\le i\le 4$, $0\le j\le 5$.
\item[(4)] $d_i d_j$ with $i\in \{1,5\}$, $j\in \{0,2,4\}$.
\item[(5)] $d_i d_i d_j$ with $i\in \{1,5\}$, $j\in \{1,3,5\}$.
\end{itemize}

Similarly for $D_4$, we take as generators $z_1 = d c^{-1}$, $z_2=c^{-1}
b$, which satisfy relations
$$
\begin{aligned}
c^2-g_1 c + 1 & =0\\
z_1 c z_1' &= c^{-1}\\
z_2' c z_2 &= c^{-1}\\
Q z'_2 z'_1+z_1 c^2 z_2 &= e_1 c^{-1}\\
c z_2 c z_1 c + Q z_1' c z_2' &= e_1,
\end{aligned}
$$
where
$$
\begin{aligned}
z'_1 &:= h_1 c^{-1}-z_1\\
z'_2 &:= f_1 c^{-1}-z_2;
\end{aligned}
$$
again, this is a Gr\"obner basis for a suitable term order, compatible with
the 112 filtration.

\subsection{Proof of Theorem \ref{pbw1}}\label{pbws2}
The efficient presentations given above are each compatible with the
corresponding $F_{123}^\bullet$ filtration (except for $D_4$, for which the
filtration is $F_{112}^\bullet$); this has the immediate consequence that
Theorem \ref{pbw1} holds for those filtrations, since we have given a
presentation with relations integral over $\Bbb C[\Bbb T]$ such that the
corresponding Poincar\'e series uniformly agrees with the infinite group
case.

It remains to consider the filtrations $F_{111}^\bullet$ for $D_4$,
$F_{111}^\bullet$ and $F_{112}^\bullet$ for $E_6$ and $F_{112}^\bullet$ for
$E_7$.  In these cases, we have been unable to find presentations with
compatible Gr\"obner bases, and must therefore use more ad hoc arguments.

For $F_{112}^\bullet(E_6)$, we may take generators $c$, $d=b c^{-1}$,
$e=c^{-1}b^{-1}c^{-1}$, with $\deg(c)=0$, $\deg(d)=\deg(e)=1$.  Using the
efficient presentation for $E_6$, we can easily solve for the relations
these satisfy in degree 2; using these relations (which again have integral
coefficients), we find that the monomials
\[
c^i, c^i e^j c^m, c^i d^l c^m, c^i e^j c^k d^l c^m
\]
with $i,m\in \{0,1,2\}$, $k\in \{0,2\}$, $j,l\ge 1$, span the algebra, and
thus must form a basis (since they do so in the infinite group case, and
thus generically).  The corresponding case of Theorem \ref{pbw1} follows.

Similarly, for $F_{111}^\bullet(D_4)$, we take generators $c$,
$w_1=c^{-1}d^{-1}$, $w_2=d b$, $w_3 = b^{-1}c^{-1}$, and find that the
monomials
\[
c^i,
c^i w_1^j c^l,
c^i w_2^j c^l,
c^i w_3^j c^l,
c^i w_1^j c w_2^k c^l,
c^i w_1^j c w_3^k c^l,
c^i w_2^j c w_3^k c^l
\]
with $i,l\in \{0,1\}$, $j,k\ge 1$, span the algebra, so again must form a
basis, giving that case of Theorem \ref{pbw1}.

For $F_{111}^\bullet(E_6)$ and $F_{112}^\bullet(E_7)$, the proof is
similar; in each case we can order the monomials at degree 2 in such a way
as to obtain monic relations with coefficients in $\Bbb C[\Bbb T]$.  Since
degree 2 relations suffice for the infinite group cases, they suffice
generically, and thus suffice in general.

\subsection{Proof of Theorem \ref{genrel}}\label{assgr}
It is straightforward to use our efficient presentations to
obtain presentations of the associated graded algebra
$H_0^{123}(t,q)$: simply remove low-order terms from the
relations.  These can be simplified considerably by removing redundant
relations; we thus obtain the following presentations.  Note that for $E_6$
and $E_7$ we follow the lead of $E_8$ in replacing the generator $b$ by a
generator $d=b c^{-1}$ (i.e., by an appropriate translation in the infinite
group).

For $E_6$:
\[
\langle
c,d\mid
c^3-g_1 c^2 + g_2 c-1,
dcd - q^{-\ell} c d c^2 d c,dcdcd
\rangle
\]

For $E_7$:
\[
\langle c,d\mid
c^4 - g_1 c^3 + g_2 c^2 - g_3 c + 1,
d c^2 d,
d c d + q^{-\ell} c d c^3 d c
\rangle
\]

For $E_8$:
\[
\langle c,d\mid
c^6 - g_1 c^5 + g_2 c^4 - g_3 c^3 + g_4 c^2-g_5 c + 1,
d c^2 d,d c^3 d,d c^4 d,
d c d + q^{-\ell} c d c^5 d c
\rangle
\]

Since $dcdcd=0$ in the $E_7$ and $E_8$ algebras, we thus obtain the
claimed uniform presentation for the algebras $H_0^{123}(t,q)$.
The theorem is proved. 

{\bf Remark.} 
Similar considerations apply for the case $F_{112}^\bullet$ for
$E_6$. The associated graded algebra of $H_0^{112}(t,q)$ 
in this case has presentation
$$
\begin{aligned}
\langle
c,d\mid{}&
  c^3-g_1c^2+g_2c-1=0,
  dccd=dcd=dcce=ecd=ece=ecce=0,\\
&
  de=Q ed+g_2 dce, dce = Q^{-1} ceccdc
\rangle.
\end{aligned}
$$

\section{The surface for $q=1$}

We can use our efficient presentations for the $123$ filtrations of $E_6$,
$E_7$, $E_8$ to give explicit equations for the surfaces $S(t,1)$ with
coordinate ring $Z(t,1)$.  Since the associated graded algebra of $Z(t,1)$
is just the center of the associated graded algebra of $H(t,1)$, we can
write down the leading terms of the generators of the center, and then
simply solve for the coefficients of the lower degree terms.  We can then
solve for the resulting equation, and put it into canonical form as in
Theorem \ref{thm:Ztqgen}.  The generators $x$, $y$, $z$ that we obtain from
this process are unfortunately rather complicated, even in the simplest
case $E_6$; it turns out, however, that the equation itself can be
described much more simply.  We consider each case in turn.

For $E_6$, we obtain the equation
\[
xyz = x^3 + a_1 x^2 + a_3 x
          + y^3 + a_2 y^2 + a_4 y
                       + z^2 + a_5 z
                              + a_6,
\]
where $a_i$ are polynomials in the parameters $u_{kj}$.  The simplest way
to specify these polynomials is as follows.  Each of the triples
$u_{k1},u_{k2},u_{k3}$ multiplies to 1, and thus specifies a point on the
maximal torus ${\bold T}(SL_3)$; as the equation is invariant under
permutations of each triple, the coefficients $a_i$ are actually (virtual)
characters of $SL_3^3$.  In fact, it turns out that these virtual
characters factor through the natural map from $SL_3^3$ to the
simply-connected group $E_6$ (mapping $SL_3^3$ to a locally isomorphic
semisimple subgroup of $E_6$); equivalently, they are functions on the maximal torus
${\bold T}(E_6)$ invariant under the action of the Weyl group.  In
particular, we can express them in terms of the fundamental characters of
$E_6$, and thus obtain the following equation.
$$
\begin{aligned}
x y z
&=
 x^3+(\chi_1) x^2+(\chi_3-\chi_2) x
+y^3+(\chi_2) y^2+(\chi_4-\chi_1) y\\
&\phantom{=}+z^2+(\chi_5-6) z+(\chi_6-3 \chi_5+9),
\end{aligned}
$$
where $\chi_1$ through $\chi_6$ are the fundamental characters.  We note
that each coefficient has a different fundamental character as its leading
term, which is associated to the coefficient as follows.  We assign a
simple root of $E_6$ to each monomial with nonconstant coefficient in such
a way that each of the three polynomials corresponds to a different leg of
the Coxeter diagram, in order of increasing degree.  In particular, the
constant term corresponds to the central root, the coefficients of degree 1
in x,y,z correspond to the roots adjacent to the center, and so forth.
(The labelling we have chosen for the roots can thus be read off from the
equation.)  Note that although the coefficients above are characters on the
simply connected group, the isomorphism class of the surface depends only
on the image in the adjoint group; the center of $E_6$ acts by
$(x,y,z)\mapsto (\zeta_3 x,\zeta_3^{-1}y,z)$.  Thus the surfaces $S(t,1)$
are parametrized by the quotient by the Weyl group of the adjoint torus;
this agrees with Theorem 8.4 of \cite{Lo}.

As a corollary, we find that the map from ${\bold T}(E_6)$ to the canonical
equation is Galois with Galois group $W(E_6)$; similarly the Galois group
of the normal closure of the map from ${\bold T}(SL_3)^3$ is the semidirect
product of $W(E_6)$ by an elementary abelian group $(1/3)Q/P$ of order
$3^5$ (the quotient by the weight lattice $P$ of $1/3$ times the root
lattice $Q$, the latter being the closure of the weight lattice of $A_2^3$
under $W(E_6)$).

Similarly, for $E_7$, the parameters $u_{kj}$ naturally specify elements of
${\bold T}(SL_2\times SL_4^2)$, and the coefficients of the equation,
virtual characters of $SL_2\times SL_4^2$, factor through the natural map
to the simply-connected group $E_7$.
We obtain the equation
$$
\begin{aligned}
xyz 
={}&
 x^4+\chi_1
 x^3+(\chi_2-2\chi_3+23)x^2+(\chi_4-\chi_1\chi_3-5\chi_6+29\chi_1)x\\
 {}+{}&y^3+(\chi_3-25)y^2+(\chi_5-\chi_2-16\chi_3+206)y\\
 {}+{}&z^2+(\chi_6-6\chi_1)z\\
 {}+{}&\chi_7+\chi_3^2-3\chi_6\chi_1+9\chi_1^2-10\chi_5+9\chi_2+62\chi_3-558,
\end{aligned}
$$ where the $\chi_i$ are fundamental characters of the simply connected
group of type $E_7$, associated to the corresponding terms in the same way
as for $E_6$.  The map from ${\bold T}(E_7)$ to the canonical equation is
Galois with Galois group $W(E_7)$; the map from ${\bold T}(SL_2\times
SL_4^2)$ has normal closure with Galois group $W(E_7)\ltimes
[(1/4)Q/P]$, where $P$ is the weight lattice of $E_7$ and $Q$
is the root lattice.

Finally, for $E_8$, the algebra is parametrized by ${\bold T}(SL_2\times
SL_3\times SL_6)$, but the center is parametrized by ${\bold T}(E_8)$; we
obtain the equation
\[
xyz
=
x^5 + a_1 x^4 + a_3 x^3 + a_5 x^2 + a_7 x +
y^3 + a_2 y^2 + a_6 y +
z^2 + a_4 z
+ a_8,
\]
where
$$\begin{aligned}
a_1=\chi_1&{}-248;\\
a_2=\chi_2&{}-25\chi_1+2325;\\
a_3=\chi_3&{}-3\chi_2-170\chi_1+23405;\\
a_4=\chi_4&{}-6\chi_3-35\chi_2+920\chi_1-57505;\\
a_5=\chi_5&{}-7\chi_4-135\chi_3-2\chi_2\chi_1+580\chi_2+23\chi_1^2+7652\chi_1-955978;\\
a_6=\chi_6&{}-\chi_5-28\chi_4+170\chi_3-16\chi_2\chi_1+2006\chi_2\\
&{}+206\chi_1^2-51436\chi_1+2401694;\\
a_7=\chi_7&{}-13\chi_6-104\chi_5-5\chi_4\chi_1+1045\chi_4-\chi_3\chi_2+29\chi_3\chi_1+4145\chi_3\\
&{}+2\chi_2^2+359\chi_2\chi_1-45708\chi_2-4444\chi_1^2+275989\chi_1+4532634;\\
a_8=\chi_8&{}-58\chi_7-10\chi_6\chi_1+1245\chi_6+9\chi_5\chi_1+2177\chi_5-3\chi_4\chi_3\\
&{}-17\chi_4\chi_2+741\chi_4\chi_1-65323\chi_4+9\chi_3^2+161\chi_3\chi_2-4405\chi_3\chi_1\\
&{}+189168\chi_3+\chi_2^2\chi_1+192\chi_2^2+62\chi_2\chi_1^2-38134\chi_2\chi_1+2537119\chi_2\\
&{}-558\chi_1^3+494091\chi_1^2-52476655\chi_1+1484285983.
\end{aligned}$$
The Galois groups are $W(E_8)$ from ${\bold T}(E_8)$ and
$W(E_8)\ltimes [(1/6)Q/P]$ (where $P=Q$ is the weight and root
lattice of $E_8$) from
${\bold T}(SL_2\times SL_3\times SL_6)$.

Note that a similar expression holds for the $D_4$ case.

{\bf Remark.} The fact that the Galois groups 
of coverings defined above are Weyl groups of $E_6,E_7,E_8$ also
follows (without computation) 
from the results of Section 7 and the Arnold-Brieskorn 
theorem, saying that the monodromy group of a simple singularity 
is the corresponding Weyl group. 

\medskip
In addition to suggesting that the canonical equations should have a
group-theoretical interpretation, the above form for the equations has a
particularly interesting consequence.  One of the most important structures
on a del Pezzo surface is the collections of lines on the surface (e.g.,
the 27 lines on a cubic surface).  Normally, these lines are only defined
over an extension field; in our case, however, it turns out that every line
on the surface is actually rational over ${\bold T}(E_k)$, so in particular
is rational in the roots of the three minimal polynomials.  Each line
intersects the curve at infinity; since the smooth part of the curve at
infinity has a natural multiplicative group structure (more precisely, a
natural divisor class of degree 3,2,1 for $E_6$, $E_7$, $E_8$), we obtain
an element of this group for each line (up to global inversion and, for
$E_6$, multiplication by a 3rd root of unity (respectively multiplication
by $-1$ for $E_7$)), each of which is a Laurent monomial in the roots of
the minimal polynomials.

For $E_6$, we obtain 27 lines, corresponding to the ratios (shortest
weights of $E_6$)
\[
b_i/a_j, c_i/b_j, a_i/c_j;\quad 1\le i,j\le 3,
\]
where $a_i$, $b_i$, $c_i$ are the roots of the minimal polynomials of $a$,
$b$, $c$, respectively.

For $E_7$, we obtain 56 lines, corresponding to the ratios (shortest
weights of $E_7$)
$$\begin{aligned}
a_i/b_jb_k,a_i/c_jc_k&;1\le i\le 2, 1\le j<k\le 4\\
b_i/c_j,c_i/b_j&; 1\le i,j\le 4.\\
\end{aligned}$$

Finally, for $E_8$, we obtain 240 lines, corresponding to the ratios (roots
of $E_8$)
$$\begin{aligned}
a_i/a_j&;1\le i,j\le 2,i\ne j\\
b_i/b_j&;1\le i,j\le 3,i\ne j\\
c_i/c_j&;1\le i,j\le 6,i\ne j\\
a_i/c_jc_kc_l&;1\le i\le 2,1\le j<k<l\le 6\\
b_i/c_jc_k,c_jc_k/b_i&;1\le i\le 3,1\le j<k\le 6\\
a_ib_jc_k,1/a_ib_jc_k&;1\le i\le 2,1\le j\le 3,1\le k\le 6.
\end{aligned}$$

Furthermore, we find that the surface is singular if and only if one or
more of the roots of $E_k$ vanishes.  The most singular case corresponds to
the identity element of the torus, in which case $a$, $b$, and $c$ are all
required to be unipotent; we readily verify in each case that the surface
has a singularity of type $E_k$ in this case.  Note in particular that for
$E_8$ the surface for the unipotent case has equation
\[
xyz = x^5+y^3+z^2;
\]
all of the lower degree terms vanish.

We conjecture that when $Q$ is a root of unity of order $k$ prime to $l$,
the coefficients of the corresponding canonical equation are simply given
by composing the above functions with the $k$-th power map on the original
torus; we have checked this for $E_6$, $k=2$, as well as the corresponding
statement for $D_4$, $k=3$.

Any $\ell$-dimensional representation of the algebra corresponds to a point
on the corresponding surface; this thus gives rise to two natural (open)
questions.  First, which representations correspond to singular points?
(Reducibility appears to be sufficient, but is not necessary.)  Second,
which representations correspond to points on a line of the surface?  For
the latter question, it appears from experiments that the family of
representations associated to a line can be parametrized in such a way that
$a$, $b$, $c$, and all their powers are linear functions in the parameter,
but it is unclear why this should be so.

\section
{Appendix 1: Generalized double affine Hecke algebras 
and multiplicative preprojective algebras, 
 by W. Crawley-Boevey and P. Shaw}

Let $K$ be a field. Let $w=(w_1,\dots,w_k)$ be a collection of positive integers,
let $\mu\in K^*$, and let $\xi_{ij}\in K^*$ ($1\le i\le k$, $1\le j\le w_i$).
Following a question of Etingof, we show that the associative
$K$-algebra $A_{w,\mu,\xi}$ with generators $x_1,\dots,x_k$ and relations
\begin{gather*}
x_1 x_2 \dots x_k = \mu 1,
\\
(x_i - \xi_{i1}1)(x_i - \xi_{i2}1)\dots (x_i - \xi_{i,w_i}1) = 0,
\quad (i=1,\dots,k),
\end{gather*}
is isomorphic to $e\Lambda^q e$ for a suitable multiplicative
preprojective algebra $\Lambda^q$ and a suitable idempotent $e$.
Note that the generalized double affine Hecke algebra $H(t,q)$
is isomorphic to $A_{w,\mu,\xi}$ for
suitable $(w,\mu,\xi)$. There is a corresponding result for
the algebra with relation $x_1+\dots+x_k=\mu 1$ in terms
of the deformed preprojective algebra, see \cite{Me}.
Note that by rescaling one of the $x_i$, and the corresponding $\xi_{ij}$,
one may assume that $\mu=1$. We make this assumption from now on.

We use the notation of \cite[\S8]{CBS},
introducing a quiver $Q_w$ with vertex set
$I = \{0\} \cup \{ [i,j] \mid 1\le i\le k, 1\le j \le w_i-1\}$,
an element $q\in (K^*)^I$, and an ordering $<$ on the arrows of $\overline{Q_w}$.
Let $\Lambda^q$ be the corresponding multiplicative preprojective algebra.
We denote by $e_v$ the idempotent corresponding to the trivial path at a vertex $v\in I$.

\begin{lemma}\label{isoo}
$A_{w,1,\xi} \cong e_0 \Lambda^q e_0$.
\end{lemma}

\begin{proof}
By \cite[Lemma 8.1]{CBS}, $e_0 \Lambda^q e_0$ is spanned by the
paths in $\overline{Q}$ with head and tail at 0. In fact it is
generated by the paths $a_{i1}a_{i1}^*$ ($i=1,\dots,k$), because any
arrows $a_{ij},a_{ij}^*$ with $j$ maximal which occur in a path, must occur
as part of a product $a_{ij}a_{ij}^*$, and if $j>1$ this
can be rewritten as $q_{i,j-1} a_{i,j-1}^* a_{i,j-1} + (q_{i,j-1}-1) e_{[i,j-1]}$.
The elements $y_i = \xi_{i1}(a_{i1}a_{i1}^* + e_0)$ clearly satisfy
$y_1 y_2 \dots y_k = e_0$, and a calculation similar to \cite[\S3]{CB1}
shows that $(y_i - \xi_{i1}e_0)(y_i - \xi_{i2}e_0)\dots (y_i - \xi_{i,w_i}e_0) = 0$.
Thus there is a surjective algebra homomorphism
$\theta:A_{w,1,\xi}\to e_0\Lambda^q e_0$ sending $x_i$ to $y_i$.
To show that $\theta$ is an isomorphism, it suffices to show that
any $A_{w,1,\xi}$-module $M$ can be obtained by restriction from
an $e_0\Lambda^q e_0$-module. Let $X$ be the representation of
$\overline{Q}$ with $X_0=M$,
\[
X_{[i,j]} = (x_i - \xi_{i1}1)(x_i - \xi_{i2}1)\dots (x_i - \xi_{ij}1)M,
\]
and with $a_{ij}$ the inclusion and $a_{ij}^*$ multiplication by $(x_i - \xi_{ij}1)/\xi_{ij}$.
Clearly $X$ defines a $\Lambda^q$-module, and $e_0 X=M$, as desired.
\end{proof}

\section{Appendix 2: Generalized double affine Hecke algebras 
and multiplicative preprojective algebras, continued}

In this appendix we will use the results of the main part of the
paper and of Appendix 1 to prove some results about the structure
of multiplicative preprojective algebras of \cite{CBS}, in
particular in the case of an affine quiver $Q$ of type $\tilde{D}_4$,
$\tilde{E}_6$, $\tilde{E}_7$,
$\tilde{E}_8$. 

\subsection{Starlike quivers}
In this subsection we give some results about 
multiplicative preprojective algebras for starlike quivers; 
they can be rather easily deduced from the paper
\cite{CBS}. 

We retain the notation of \cite{CBS}. 
In particular, given a function $q$ on the set of vertices of $Q$
with values in $\Bbb C^*$, $\Lambda^q$ denotes the corresponding
multiplicative preprojective algebra. 

\begin{proposition}\label{typeA} Let $Q$ be the quiver of type $A_m$, with
vertices labeled $1,...,m$. In this case the algebra 
$\Lambda^q$ is finite dimensional, and
it is zero unless $\prod_{p=i}^j q_p=1$ for some $i\le j$. 
\end{proposition}

\begin{proof} Let $a_i$ be the edges from $i$ to $i+1$ and
$a_i^*$ from $i+1$ to $i$. The algebra $\Lambda^q$ 
is the quotient of the path algebra of the double of $Q$ by the
relations
$$
1+a_1^*a_1=q_1,\ q_2(1+a_2^*a_2)=1+a_1a_1^*,...,
q_m=1+a_{m-1}a_{m-1}^*.
$$
This algebra is isomorphic to a usual (deformed) preprojective algebra 
for the quiver $Q$. In particular, 
it is finite dimensional (as $Q$ is of finite Dynkin type; see
\cite{CBH}). 

Now, if $\Lambda^q\ne 0$, then it must have a finite dimensional
irreducible representation. Therefore, by Theorem 1.9 in
\cite{CBS}, $q^\alpha=1$ for some positive root $\alpha$, as desired. 
\end{proof}

Let $Q$ be any starlike quiver, with vertices
labeled by pairs $(j,k)$, where $k=1,...,m$ is the number of the
leg, and $j=1,...,d_k-1$ the number of the vertex on the $k$-th
leg, enumerated from the nodal vertex. 

\begin{proposition}\label{starl} The algebra 
$B:=\Lambda^q/\Lambda^q e_0\Lambda^q$ is finite dimensional. 
Furthermore, if $\prod_{p=i}^j q_{(p,k)}\ne 1$ for any
$k=1,...,m$ and $1\le i\le
j<d_k$, then this algebra is zero and hence
we have a natural Morita equivalence between 
$\Lambda^q$ and $e_0\Lambda^q e_0$.  
\end{proposition}

\begin{proof}
The algebra $B$ is a cyclic $\Lambda^q$-module in which $e_0$ acts by
zero. Thus $B$ is supported at the non-nodal vertices of $Q$. 
Hence $B$ can be regarded as a cyclic module over the direct sum  
of the multiplicative preprojective algebras for quivers of types
$A_{d_k-1}$. By Proposition \ref{typeA}, this implies that $B$ is
finite dimensional. Also, if the condition on parameters holds, these
algebras are zero, and hence $B=0$, as desired. 
\end{proof}

Let $L^q$ be the quotient of the path algebra of the double of $Q$ by the
relations of the multiplicative preprojective algebra except 
the nodal relation. 

\begin{lemma}\label{finge} (i) Let $i$ be a vertex 
of $Q$. Then $e_i L^q e_0$ is a cyclic right module over $e_0 L^q e_0$ 
generated by the straight (=shortest) path connecting the nodal
vertex $0$ with $i$. 

(ii) The same is true for $e_i \Lambda^q e_0$ as a right module over 
$e_0 \Lambda^q e_0$. 
\end{lemma}

\begin{proof} It is clear that (ii) follows from (i), so it
suffices to prove (i). 
Let $\gamma$ be any path leading from $0$ to the vertex $i$. 
We need to show that $\gamma$ is a linear combination of straight
paths over $e_0 L^q e_0$. If $\gamma$ visits more than one leg,
then it is a product of a shorter path with an element of $e_0
L^q e_0$, so by induction (in the length of the path) 
we may assume that $\gamma$ is entirely contained in
one leg, and hence, without loss of generality, that $Q$ has only
one leg to begin with. We will assume that the orientation of $Q$
is toward the nodal vertex. If $\gamma$ is not straight, it contains a
subpath $aa^*$ for some edge $a$. If $a$ ends at the vertex
$0$, then $\gamma$ is a product of a shorter path and a path
starting and ending at $e_0$, and by induction we are done. 
Otherwise, we may use the relation at the head of $a$ 
to replace $aa^*$ with $c_1 b^*b+c_2$ where $b$ is the edge 
that begins at the head of $a$ and is directed toward the nodal
vertex, and $c_1$ and $c_2$ are constants. This allows us to 
represent $\gamma$ as a linear combination of a shorter path
and a path of the same length but smaller sum of distances of
the vertices passed to the nodal vertex. This proves that
the straight path is a generator, as desired. 
\end{proof}

\begin{corollary} \label{idem}
If the condition on $q$ in 
Proposition \ref{starl} holds, then for each $i$
$e_i \Lambda^q e_0=p_i e_0\Lambda^q e_0$, where $p_i\in
e_0\Lambda^q e_0$ are certain idempotents.
Hence $\Lambda^q=\oplus_{i,j}p_i(e_0\Lambda^q e_0)p_j$. 
\end{corollary} 

\begin{proof} Indeed, by Proposition \ref{starl} $e_i \Lambda^q
e_0$ is a projective module over $e_0\Lambda^q e_0$, 
while by Lemma \ref{finge} it is a quotient (hence a direct
summand) of a free rank 1 module. This implies the result. 
\end{proof} 

In fact, it is easy to see that if $i$ belongs to the 
$k$-th leg then the idempotent $p_i$ 
is a polynomial of the element $x_k$
defined in Appendix 1 projecting to the direct sum of 
${\rm distance}(i,0)$ eigenspaces of $x_k$
(note that the condition on $q$ in 
Proposition \ref{starl} is equivalent to saying that 
the elements $x_k$ have distinct eigenvalues). 

\subsection{Affine quivers}

Let us now apply the results of the previous subsection to affine
quivers. Assume that $Q$ is of type $\tilde{D}_4$, 
$\tilde{E}_6$, $\tilde{E}_7$,
$\tilde{E}_8$. Let $\delta$ be a basic imaginary root. 
By Lemma \ref{isoo}, in this case the algebra
$e_0\Lambda_q e_0$ is isomorphic to the generalized double affine
Hecke algebra $H=H(t,q')$ for certain $t,q'$
related to $q$ by a simple transformation (for example, $q'=q^\delta$).

\begin{corollary}\label{mori} 

(i) The Gelfand-Kirillov dimension of $\Lambda^q$ is $2$. 

(ii) If the condition on $q$ in Proposition \ref{starl} holds, then 
the algebra $\Lambda^q$ is naturally Morita equivalent to $H$. In
this case, $\Lambda_q=\oplus_{i,j} p_iHp_j$.  
for the corresponding $t$, $q'$. 

(iii) Suppose that in the situation of (ii), $q'$ is a root of unity such that 
$(q')^\ell$ has order $N$, and $t$ is so generic that $H$ is an Azumaya
algebra (of degree $\ell N$; see Section 5). Then $\Lambda^q$ is 
an Azumaya algebra of degree $hN$, where $h$ is the Coxeter number of the
corresponding Dynkin diagram. 
\end{corollary}

\begin{proof} (i) Let us put a filtration on $\Lambda^q$,
which assigns degree $0$ to vertex idempotents and degree $1$ to
edges. It is easy to see that this filtration extends 
the length filtration on $H$ defined in Section 5. 
This together with Proposition \ref{finge} implies that 
the algebra $\Lambda^q$ exhibits the quadratic growth 
in this filtration, hence the result. 

(ii) This follows from Corollary \ref{idem}.

(iii) It is easy to see that 
for a given $k$ the multiplicities of all eigenvalues of
the element $T_k\in H$ in an irreducible representation of $H$
are the same, and equal $\ell N/d_k$. Therefore, 
the rank of the idempotent $p_i$ in such a representation 
is equal to $N$ times the $i$-th coordinate $\delta$. 
The sum of such coordinates is the
Coxeter number. Thus the statement follows from (ii). 
\end{proof}

\end{document}